\documentclass[12pt]{amsart}
\usepackage[utf8]{inputenc}
\usepackage[usenames,dvipsnames]{xcolor}
\usepackage{tikz}
\usepackage{subfigure}
\usepackage{fancyhdr}
\usepackage{txfonts}
\usepackage{epsfig}
\usepackage{mathrsfs}
\usepackage{amssymb}
\usepackage{latexsym}
\usepackage{amsmath}
\usepackage{amsfonts}
\usepackage{amsbsy}
\usepackage{amscd}
\usepackage{amsthm}
\usepackage{verbatim}
\usepackage{color,tikz}
\usetikzlibrary{calc}
\usetikzlibrary{decorations.pathreplacing}
\usetikzlibrary{decorations.markings}
\usetikzlibrary{intersections}

\usepackage{amsbsy,calc}
\usepackage{ifthen,cite}
\usepackage{appendix}
\usepackage{latexsym,amsthm,amssymb,amsmath,amsfonts,graphicx,epstopdf,subfigure,float}
\usepackage{bbm}
\usepackage{pifont}
\usepackage{tikz,xcolor}

\newtheorem{theorem}{Theorem}[section]
\newtheorem{thm}[theorem]{Theorem}
\newtheorem{pro}{Proposition}[section]
\newtheorem{cor}[pro]{Corollary}
\newtheorem{lemma}[pro]{Lemma}
\newtheorem{remark}[pro]{Remark}
\newtheorem{defi}{Definition}[section]

\newtheorem{example}{Example}[section]
\numberwithin{equation}{section}

\def\b{\bigskip}

\def\cal{\mathcal }

\def\mathscr{\mathcal }

\def\be{\mathbf e}
\def\ba{\mathbf a}

\def\bd{\mathbf d}
\def\bs{\mathbf s}

\def\b0{\mathbf 0}

\newcommand{\bx}{{\mathbf{x}}}
\newcommand{\by}{{\mathbf{y}}}
\newcommand{\bz}{{\mathbf{z}}}
\newcommand{\bu}{{\mathbf{u}}}
\newcommand{\bv}{{\mathbf{v}}}
\newcommand{\bb}{{\mathbf{b}}}

\newcommand{\SD}{{\mathcal{D}}}

\newif\ifdraft
\drafttrue
\ifdraft
\else
\fi
\begin{document}

\title{Topology automaton and H\"older equivalence of Bara\'nski carpets}

\author{Yunjie Zhu}
\address{Department of Mathematics and Physics, Hubei Polytechnic University, Huangshi,435000, China.}
\email{yjzhu\_ccnu@sina.com}

\author{Liang-yi Huang$^\dag$}
\address{College of Science, Wuhan University of Science and Technology, Wuhan, 430070, China}
\email{liangyihuang@wust.edu.cn}

\author{chunbo Cheng}
\address{Department of Mathematics and Physics, Hubei Polytechnic University, Huangshi,435000, China.}
\email{bccheng@hbpu.edu.cn}

\date{\today}
\thanks{This work is supported by NSFC No. 12401110 and No. 62206088.}
\thanks{$\dagger$ Corresponding author.}
\maketitle{ }

\begin{abstract}
	The study of Lipschitz equivalence of fractals is a very active topic in recent years. In 2023, Huang \emph{et al.} (\textit{Topology automaton of self-similar sets and its applications to metrical classifications}, Nonlinearity \textbf{36} (2023), 2541-2566.) studied the H\"older and Lipschitz equivalence
of a class of p.c.f. self-similar sets which are not
 totally disconnected. The main tool they used is the so called topology automaton.
	In this paper,  we define topology automaton for Bara\'nski carpets,
	 and we show that the method used in Huang \emph{et al.} still works for the self-affine and non-p.c.f. settings.
	As an application, we obtain a very general sufficient condition for Bara\'nski carpets to be H\"older (or Lipschitz) equivalent.
\end{abstract}

\section{Introduction}

To determine whether two  metric spaces are homeomorphic, quasi-symmetric or Lipschitz equivalent is an important
and active topic
in analysis.
In recent years, there are a lot of works devoted to the  quasi-symmetric equivalence of fractal sets, see \cite{BT, Bonk, Solomyak10, TW06, Kigami_2014, Dang_2021}.
Since 1990, the study of   bi-Lipschitz classification of self-similar sets becomes hot and abundant results have been obtained, see \cite{FM92, DS, FanRZ15,LuoL13,RRX06,RuanWX14,RaoZ15,XiXi20}. However,  most of the studies in literature focus on self-similar sets which are totally disconnected.

Lately, Huang,Wen,Yang and Zhu\cite{HWYZ23}  introduced a notion of topology  automaton for posted-critically finite (p.c.f.) fractals. Using this new tool,  they  give sufficient conditions
 for H\"older or Lipschitz equivalence   of  a class of self-similar sets, called fractal gaskets,   which  are not totally disconnected.

In this paper, we   define topology automaton for Bara\'nski carpets.
Then following the apporach of Huang \emph{et al.}, we use topology automaton to
study the H\"older   and Lipschitz equivalence of Bara\'nski carpets.

		Two metric spaces (or pseudo-metric spaces) $(X,d_X)$ and $(Y,d_Y)$ are said to be  \emph{H$\ddot{o}$lder equivalent}, denoted by $X\overset{\text{H\"older}}{\simeq}Y$,
		if there is a bijection $f:~X\to Y$, and  constants $s, C>0$ such that
\begin{equation}\label{eq:Holder}
 C^{-1}d_X(x_1,x_2)^{1/s} \leq d_Y \big( f(x_1),f(x_2)\big ) \leq C d_X(x_1,x_2)^s,~~~~\forall~~x_1,x_2\in X.
 \end{equation}
		In this case, we say $f$ is a \emph{bi-H$\ddot{o}$lder map} with index $s$.
		If $s=1$, we say $X$ and $Y$ are  \emph{Lipschitz equivalent}, denoted by $X\simeq Y$, and call $f$ a \emph{bi-Lipschitz map}.

First, let us recall the definition of Bara\'{n}ski carpet.
Let $n,m\geq 2$ be two integers.
Let $(f_i)_{i=0}^{n-1}$ and $(g_j)_{j=0}^{m-1}$ be two collections of contracting similarities of $[0,1]$ with positive contraction ratios such that
$$[0,1)=f_{0}[0,1)\cup \cdots \cup f_{n-1}[0,1) \text{ and } [0,1)=g_{0}[0,1)\cup \cdots \cup g_{m-1}[0,1)$$
be two partitions of $[0,1)$ from left to right. We shall call $(f_i)_{i=0}^{n-1}$ and $(g_j)_{j=0}^{m-1}$
the \emph{base horizontal IFS} and the \emph{base vertical IFS}, respectively.

Let  $\SD=\{\bd_1,\dots, \bd_N\}\subset \{0,\dots, n-1\}\times\{0,\dots, m-1\}$. We call $\SD$ a \emph{digit set}.
For each $\bd_j=(d_{j,1},d_{j,2})\in \SD$,
define
\begin{equation}\label{eq-IFS-B}
	{\varphi_j}\left({x,y}\right)= ( f_{d_{j,1}}(x), g_{d_{j,2}}(y)).
\end{equation}
Then $\varphi=\left({\varphi_j}\right)_{j=1}^N$  is a self-affine IFS, and we call its attractor $K=K_\varphi$
a \emph{Bara\'nski carpet}, see \cite{Baranski07}.

For simplicity, we denote by ${\mathcal B}(n,m, \SD)$ the collection of Bara\'nski carpets with division numbers
$n, m$ and digit set $\SD$.
If all maps in $(f_i)_{i=0}^{n-1}$ and $(g_j)_{j=0}^{m-1}$ have contraction ratios $1/n$ and $1/m$ respectively, then $K_\varphi=K(n,m,\SD)$ is called  a \emph{Bedford-McMullen carpet}, see \cite{Bedford84, McMullen84}.
If in addition $n=m$,  then $K_\varphi=K(n,\SD)$ is called  a \emph{fractal  square}, see \cite{LLR13, XiXi10}.


 Let $K\in {\mathcal B}(n,m,\SD)$ and denote  $\Sigma=\{1,\dots,N\}$.
 For $1\leq j\leq N$, we denote $K_j=\varphi_j(K)$ and we call it a first order \emph{cylinder}.
 If $\bd_j=(p,q)$, then we say $K_j$ locates in the $p$-th row and the $q$-th column.

In this paper, we confine ourselves to  Bara\'nski carpets satisfying the following  separation conditions.

We say $K$  satisfies the \emph{cross intersection condition}, if for any $i,j\in\Sigma$, $K_i\cap K_j\neq \emptyset$ implies that  $K_i$ and $K_j$ either locate in a same row, or   in a same column.
Especially,
we say   $K$ satisfies the \emph{vertical separation condition}, if  $K_i\cap K_j\neq \emptyset$  implies that $K_i$ and $K_j$  locate in a same row.

We say   $K$ satisfies \emph{top isolated condition}, if the top row of $\SD$ has only one element, say, $\SD\cap \{0,1,\dots,n-1\}\times \{m-1\}=\{\bd_{j^*}\}$,
and
$K_{j^*}\cap K_i=\emptyset$  for all $i\neq   j^*$.
(We shall call $K_{j^*}$ the \emph{top cylinder} of $K$.)


\begin{remark}
	\emph{ It is not hard to show that if $K$    satisfies the top isolated condition or the vertical separation condition, then  any non-trivial connected component of    $K$ must be a horizontal line segment.}
\end{remark}


Next, we define horizonal blocks of Bara\'nski carpets.

\begin{defi}
\emph{Let $K\in {\mathcal B}(n,m, \SD)$ be a Bara\'{n}ski carpet.}
	
\emph{(i) We call $I \subset \Sigma$ a \emph{horizontal block (H-block)} of $K$ if all $K_i$, $i\in I,$ are located in a same row,	$\bigcup_{i\in I}\varphi_i([0,1]^2)$ is connected, and $I$ is maximal with this property. We call $\# I$ the size of $I$.}
	
\emph{(ii) An H-block of size $n$ is called  a \emph{full  H-block}.}

\emph{(iii)
We call $I$ a \emph{left H-block} (resp. \emph{right H-block})  if it is not a full H-block and $\bigcup_{i\in I}\varphi_i([0,1]^2)$
intersects the left side (resp. right side) boundary of  $[0,1]^2$. }	
	
	\emph{(iv) Let  $I'$ be a left H-block and  $I''$ be a right H-block. If they  locate  in the same row, then
		we call $(I',I'')$ a \emph{H-block  pair}, and call $(\#I',\#I'')$ the size of the pair.	}
\end{defi}

	\begin{figure}[H]
		\centering
		\includegraphics[width=8.5 cm]{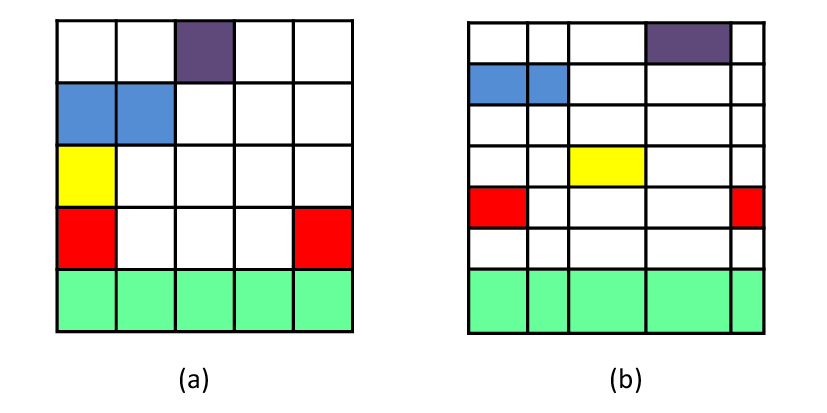}\\
		\caption{(a) $E=K(5,\cal D_E)$ is a fractal square which satisfies the top isolated condition.
 (b) $F=K(5, 7, \cal D_F)$ is a Bara\'nski carpet  which satisfies the vertical separation condition.}
		\label{fig:1}
	\end{figure}

The main result of this paper is the following.
\begin{theorem}\label{thm:main}
	Let $E\in {\mathcal B}(n, m_1,{\cal D}_E)$  and $F\in {\mathcal B}(n,  m_2,{\cal D}_F)$ be two  Bara\'{n}ski carpets satisfying
	the cross intersection condition, and assume that $E$, and also $F$, satisfies  either the top isolated condition or  the vertical separation condition.
	If there is a bijection between the H-blocks of $E$ and $F$ preserving their sizes,
	and there is a bijection between the H-block pairs of $E$ and $F$ preserving their sizes, then
	$E$ is H\"{o}lder equivalent to $F$.
	
	If in addition that  $n=m_1=m_2$, that is,  both $E$ and $F$ are  fractal squares, then $E\simeq F$.
\end{theorem}

\begin{example}\label{exam:one}
	\rm{
		Figure \ref{fig:1} illustrates two   Bara\'{n}ski carpets $E$ and $F$. It is easy to see that both $E$ and $F$ satisfy
		the cross intersection condition. Moreover, $E$ satisfies the top isolated condition and $F$ satisfies the vertical separation condition.
		Each of them contains one full H-block, one H-block with size 2, two H-blocks with size 1 and one H-block pair with size (1,1). Therefore, $E$ is  H\"{o}lder equivalent to $F$ by Theorem \ref{thm:main}. Consequently, $E$ is homeomorphic to $F$.}
\end{example}

\begin{figure}
	\centering
	
	\subfigure[ $K_1$]{
		\includegraphics[width=5cm]{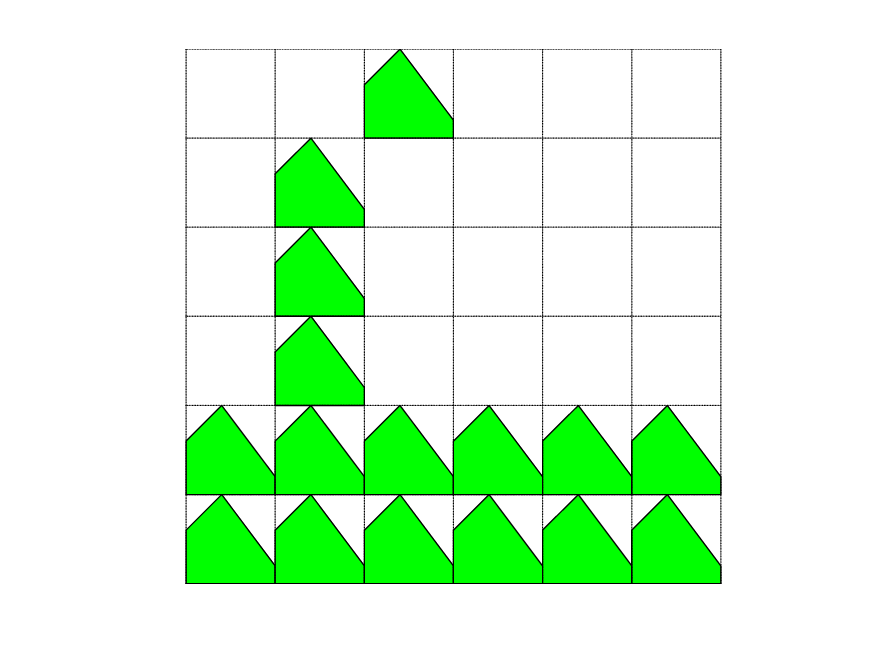}
	}
	\subfigure[ $K_2$]{
		\includegraphics[width=5cm]{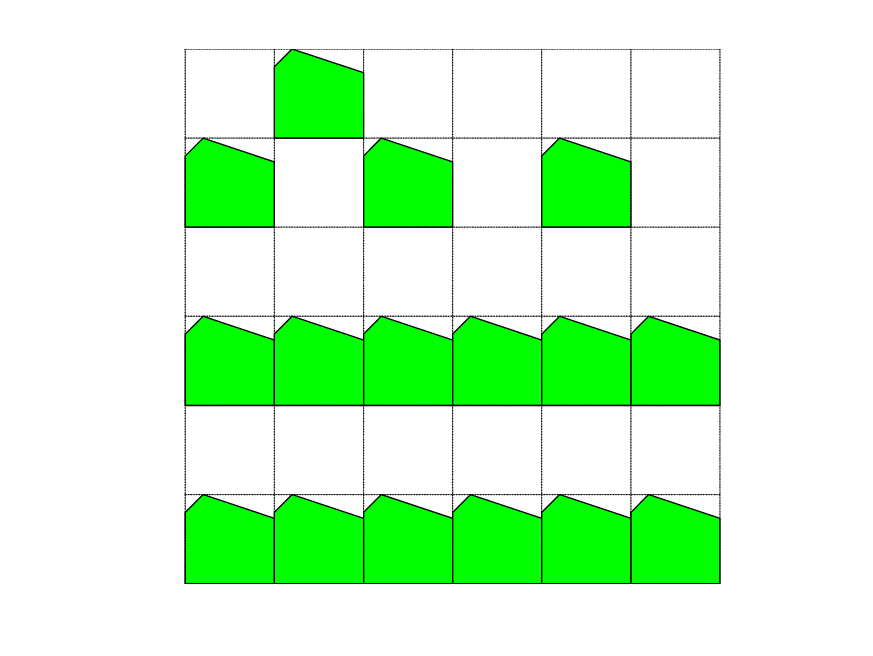}
	}\\
\subfigure[ $K_{3}$]{
		\label{fig:2:a}
		\includegraphics[width=5cm]{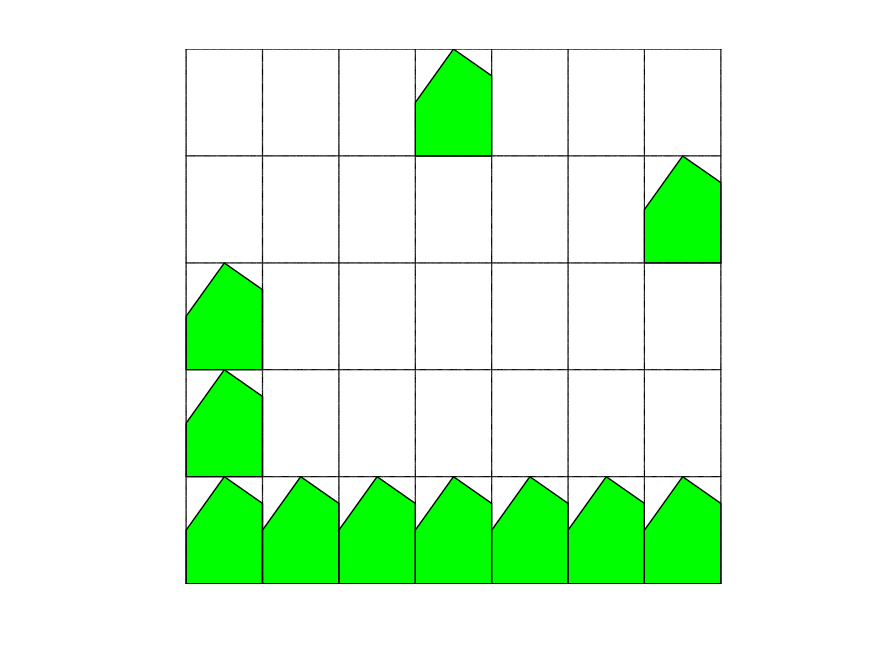}
	}
	\subfigure[$K_4$]{
		\includegraphics[width=5cm]{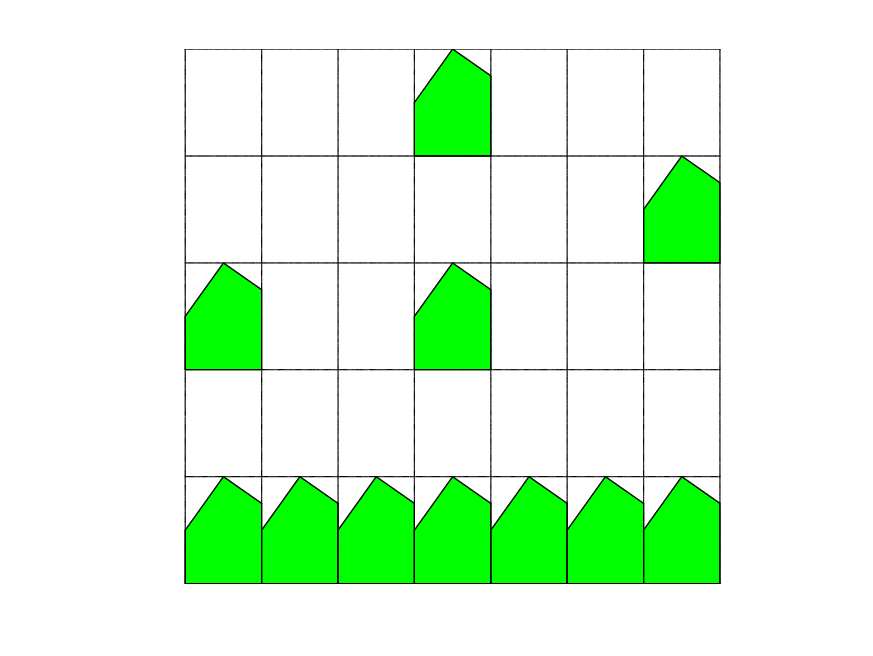}
	}\\
	\caption{$K_1\simeq K_2$ and $K_3\simeq K_4$.}\label{fig:2}
\end{figure}

\begin{example}
\emph{The essence of Theorem \ref{thm:main} is, if two cylinders in adjacent rows intersect, we can  decouple them.}

\emph{Figure \ref{fig:2} (a) and (b) illustrate two fractal squares,
where the small polygons are  convex hulls of the first order cylinders. By Theorem \ref{thm:main},  $K_1\simeq K_2$.}

\emph{Figure \ref{fig:2} (c) and (d) illustrate two Bedford-McMullen carpets.   By Theorem \ref{thm:main},
   $K_3$ and  $K_4$ are H\"older equivalent.  However, according to Rao, Xu and Zhang \cite{RXZ24},  $K_3$ and $K_4$
 are not Lipschitz equivalent, since both of them are not totally disconnected, and the fiber sequence
 of $K_3$ is not a permutation of that of $K_4$.  (Let $s_j$ be the number of cylinders of $K$ in the $j$-th row, then
 we call $(s_j)_{j=0}^{m-1}$ the \emph{fiber sequence} of $K$.)}
\end{example}

\begin{figure}[h]
	\centering
	\subfigure[ $F_{1}$]{
		\label{fig:1:a}
		\includegraphics[width=3cm]{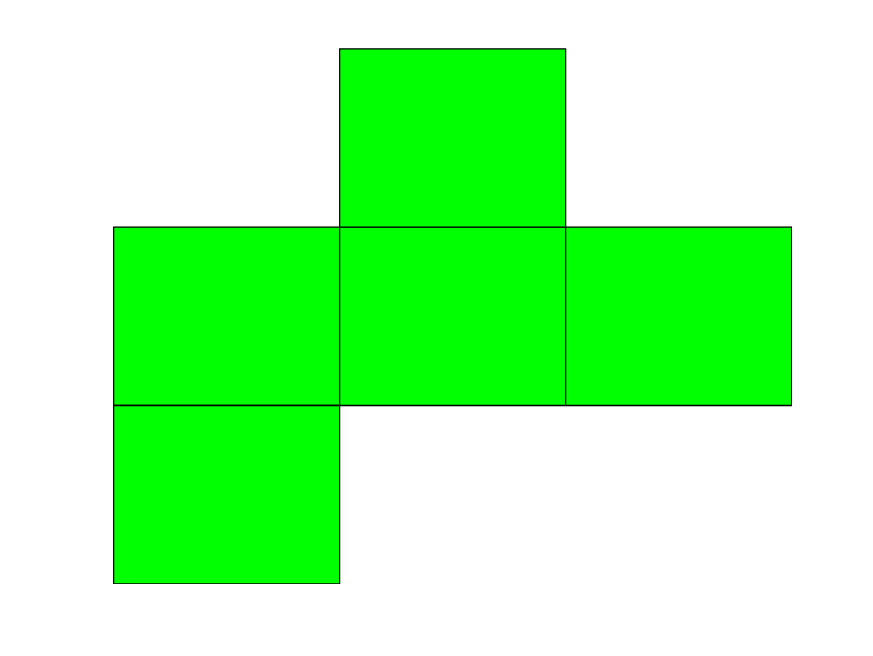}
		\includegraphics[width=3cm]{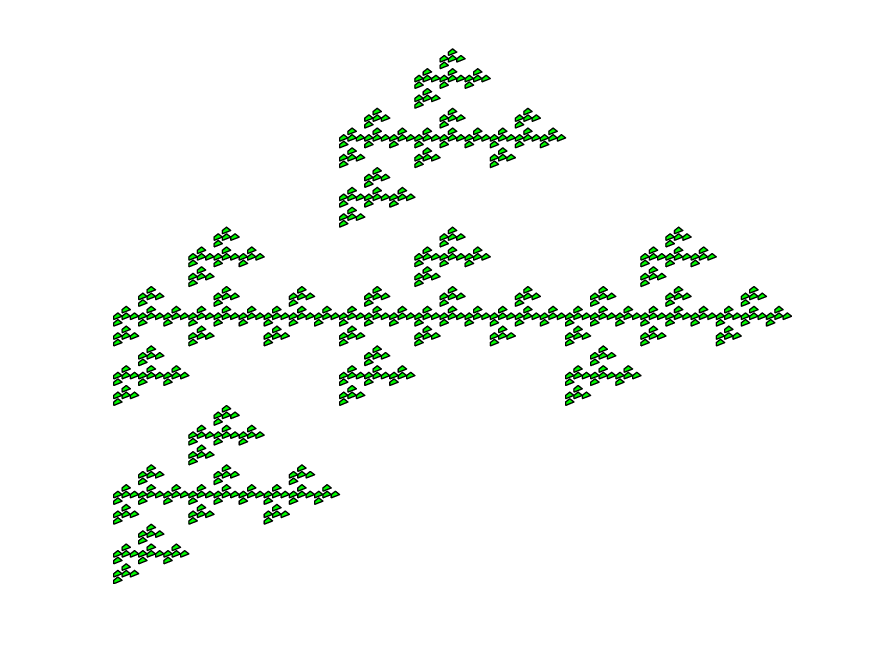}}
	\subfigure[ $F_{2}$]{
		\label{fig:1:b}
		\includegraphics[width=3cm]{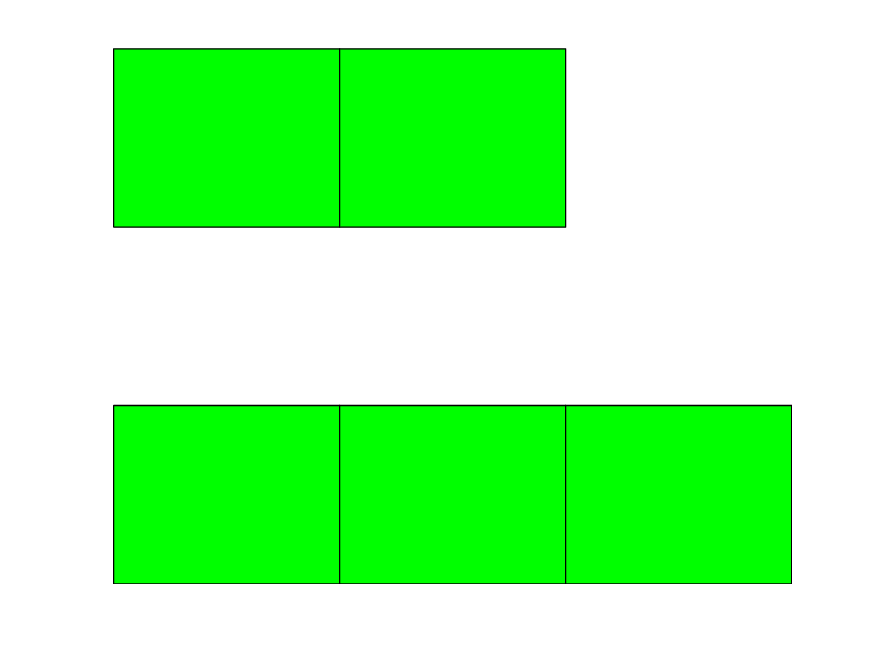}
		\includegraphics[width=3cm]{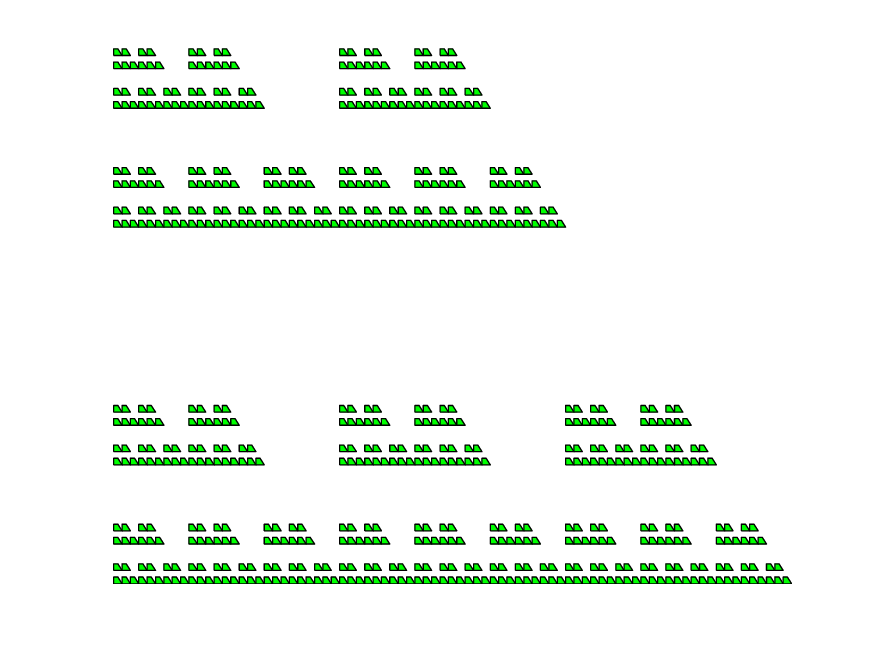}}
	\subfigure[ $F_{3}$]{
		\label{fig:1:c}
		\includegraphics[width=3cm]{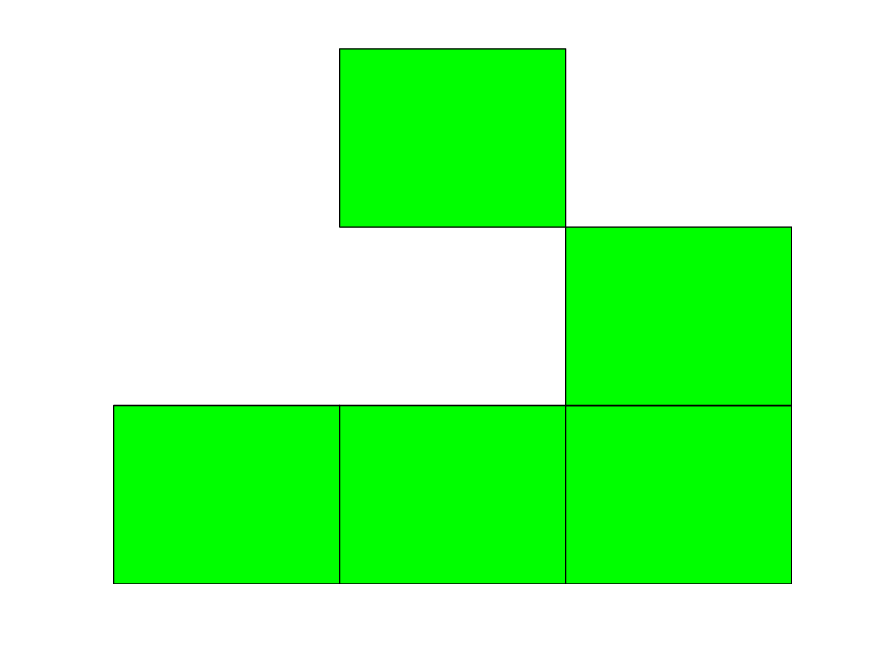}
		\includegraphics[width=3cm]{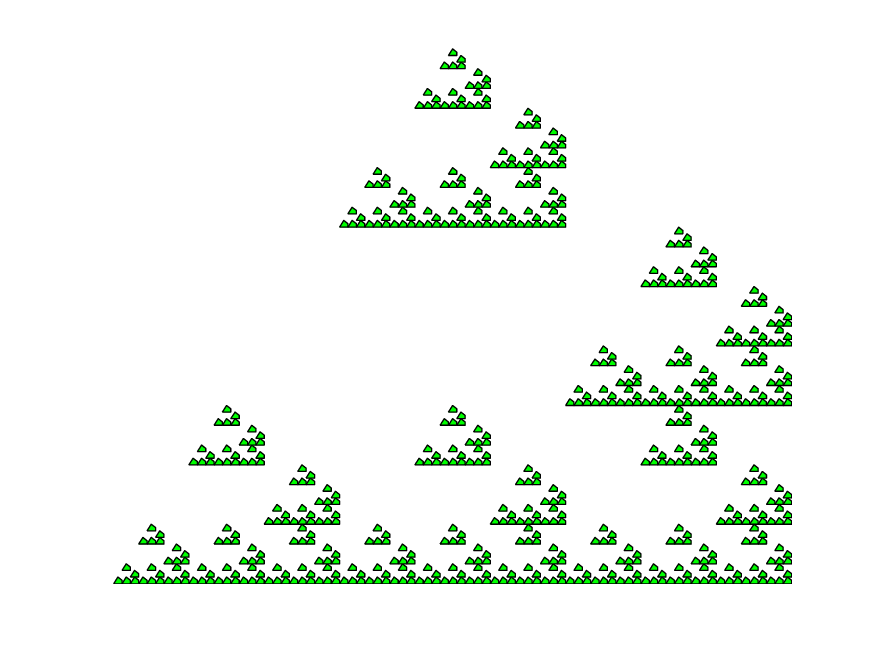}}
	\subfigure[ $F_{4}$]{
		\label{fig:1:d}
		\includegraphics[width=3cm]{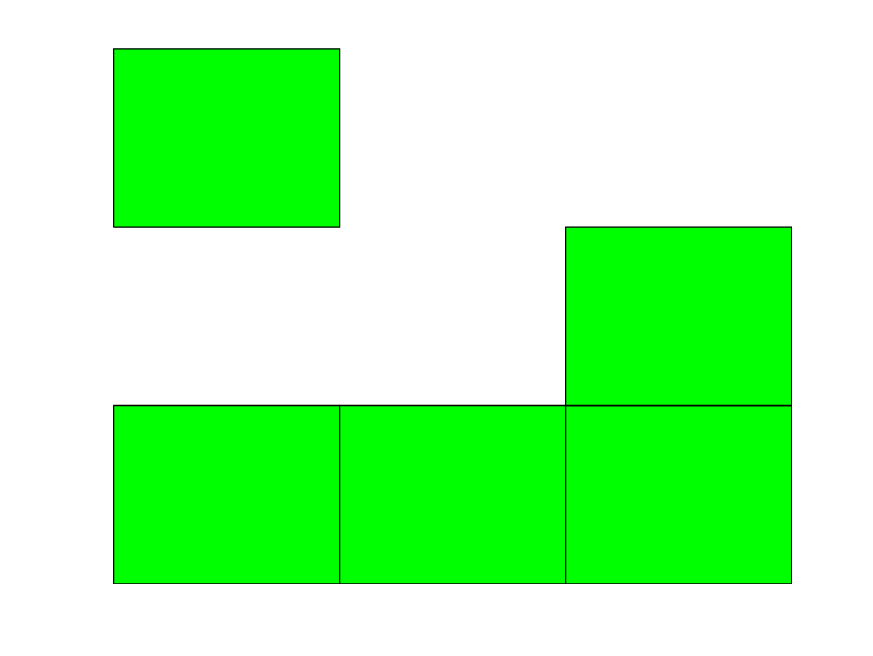}
		\includegraphics[width=3cm]{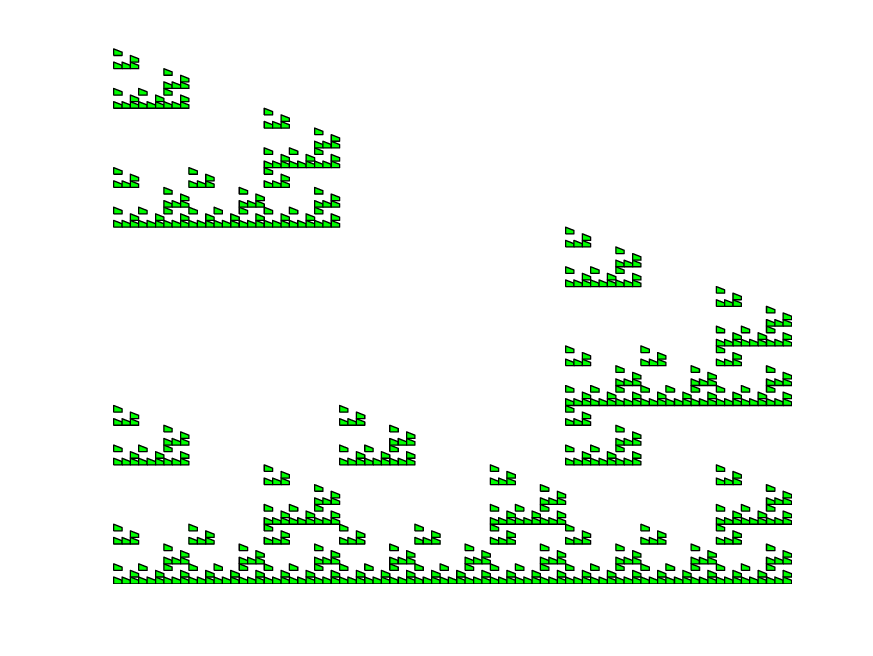}}
	
	\subfigure[ $F_{5}$]{
		\label{fig:1:e}
		
		\includegraphics[width=3cm]{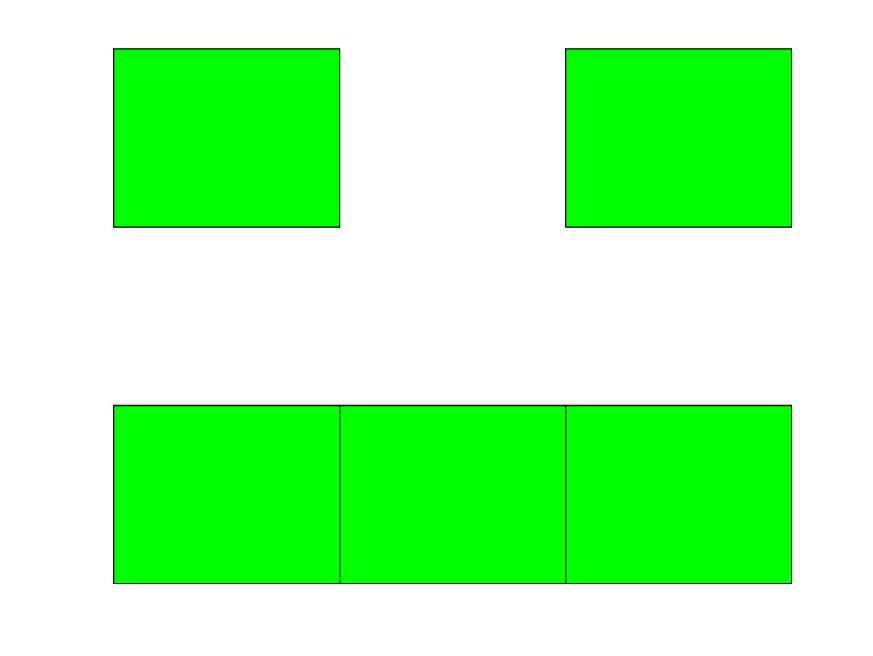}
		\includegraphics[width=3cm]{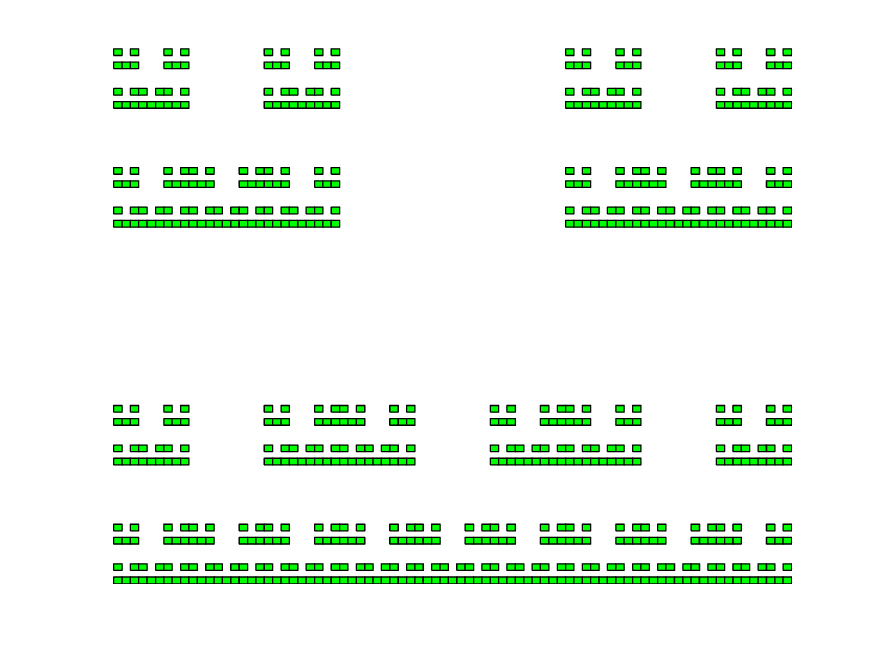}}
	\subfigure[ $F_{6}$]{
		\label{fig:1:f}
		\includegraphics[width=3cm]{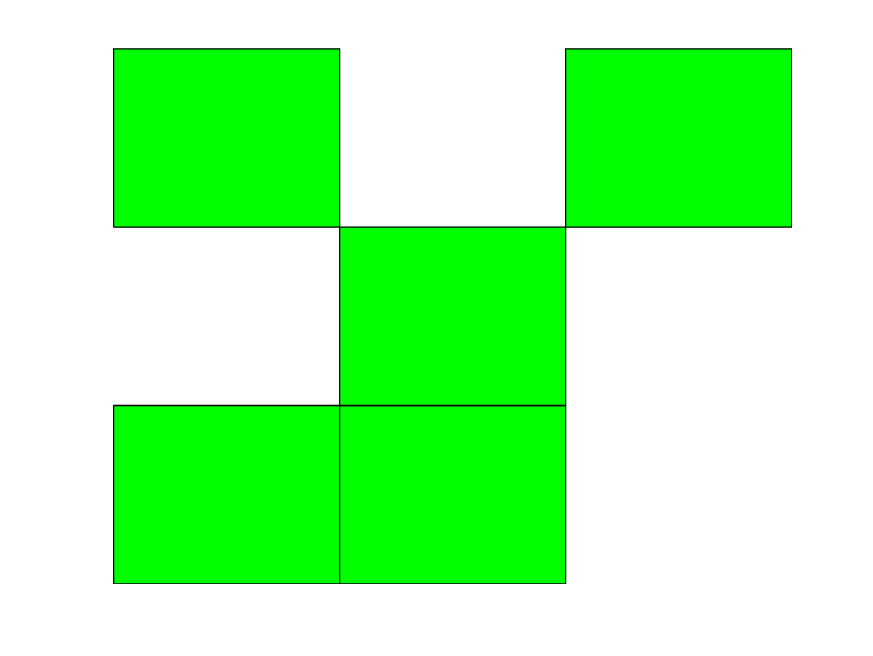}
		\includegraphics[width=3cm]{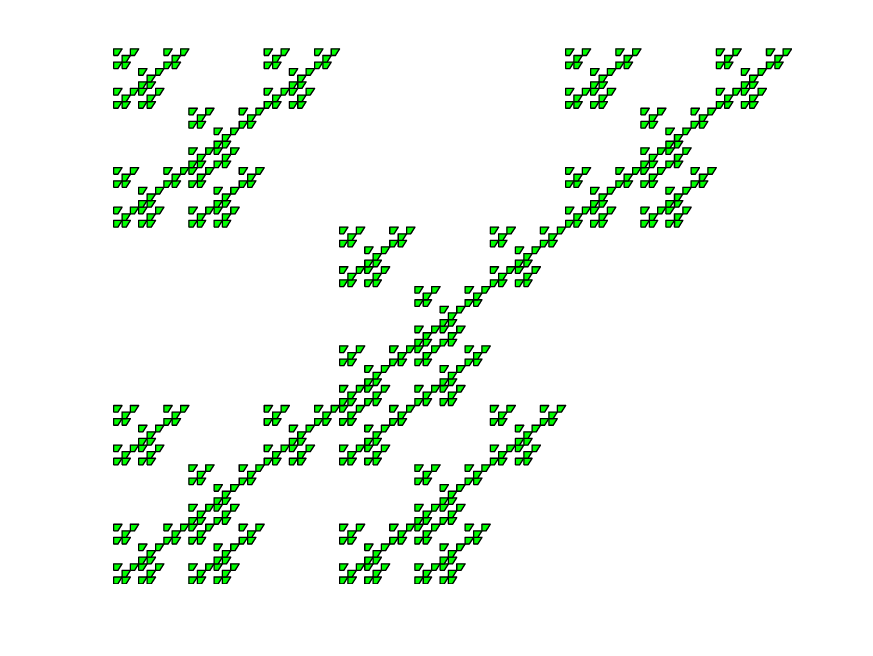}}
	
	\caption{Some fractal squares with $n=3$, $\#{\cal D}=5$.}
	\label{fig:3}
\end{figure}

\begin{example}
	[Lipschitz classification of fractal squares with expanding factor $3$]\label{exam:factor=3}
	
	\emph{Luo and Liu \cite{LuoL16} raised the question to give a complete Lipschitz classification of fractal squares with expanding factor $3$.
		Let $E=K(3, \SD_E)$ and $F=K(3, \SD_F)$ with $\#\SD_E=\#\SD_F=N$.}
	
	\emph{  The case $N\leq 4$ was settled by
		Wen \emph{et al.} \cite{WenZD12}, the case $N=6$ was settled by Rao \emph{et al.} \cite{RWZ18},
		and the cases $N=7, 8$ were settled by Ruan and Wang \cite{RuanW17}.
		Luo and Liu \cite{LuoL16} dealt with the case $N=5$, but they were not able to determine the Lipschitz equivalence relation about the 6 fractal squares depicted in Figure \ref{fig:3}. Later,
		Zhu and Rao \cite{ZhuRao21} showed that $F_1\simeq F_3$, Zhu and Yang \cite{YZ18} showed that  $F_{1}\simeq F_{2}$. }
	
	\emph{We remark that $F_1\simeq F_3$ can be obtained by Theorem \ref{thm:main}.  Clearly $F_1$ satisfies the vertical separation condition,
		and $F_3$ satisfies the top isolated condition. Moreover, both of them contain $1$ full H-block and two H-blocks of size $1$,
		therefore, $F_1\simeq F_3$.}
	
	\emph{We conjecture that $F_1\simeq F_4$, but $F_1$, $F_5$ and $F_6$ are not Lipschitz equivalent to each other.}
\end{example}

\begin{remark}
\emph{The result and proof of this paper is inspired by \cite{HWYZ23}, which dealt with a class
 of self-similar sets called	fractal gaskets.   One can easily define topology automaton for Bara\'nski carpets (see Section 3),   and the topology automaton plays a crucial role in our discussion as it does in \cite{HWYZ23}.	It is still unclear how to define topology automaton for general self-similar sets or self-affine sets.}
\end{remark}
	

	\medskip
	
	This article is organized as follows:
	In Section \ref{FSA and SN}, we introduce the pseudo-metric spaces induced by feasible $\Sigma$-automata.
We  define topology automaton of Bara\'nski carpets in Section \ref{Top-auto}, and introduce the notion of  cross automaton in Section \ref{CA}.
	In Section \ref{one Sim}, we study the one-step simplification of cross automaton.
Theorem \ref{thm:main} is proved in Section \ref{sec6}.
	In Section \ref{map g}, we recall the universal map $g$ on symbolic space given by \cite{HWYZ23}.
  Finally, we prove Theorem \ref{main} in Section \ref{pf main}, which is crucial in the proof of Theorem \ref{thm:main}.

\section{\textbf{  $\Sigma$-automaton and pseudo-metric space} }\label{FSA and SN}

In this section, we  recall some basic definitions and facts about $\Sigma$-automaton and pseudo-metric spaces,
which were introduced by Huang \emph{et al.} \cite{HWYZ23}.

\subsection{Pseudo-metric space}\label{subsec:pseudo}
~\\
\indent Let us recall the definition of pseudo-metric space, see for instance, \cite{Pepo90}.

\begin{defi}\label{def:psuedo}
\emph{
 A \emph{pseudo-quasi-metric space}
is a pair $(\mathcal{A},\rho)$ where $\mathcal{A}$ is a set and $\rho:\mathcal{A}\times\mathcal{A}\rightarrow\mathbb{R}_{\geq 0}$ satisfying for all $x,y,z\in \mathcal{A}$, it holds that  $\rho(x,x)=0$, $\rho(x,y)=\rho(y,x)$, and }
\emph{
$$\rho(x,z)\le C'\big(\rho(x,y)+\rho(y,z)\big),$$
  where $C'\ge 1$ is constant independent of $x,y,z$.}

\emph{If in addition $x\neq y$ implies $\rho(x,y)>0$, then   $({\mathcal A}, \rho)$ is called
a \emph{pseudo-metric space}. }
\end{defi}

There is a standard way to construct a pseudo-metric space from a pseudo-quasi-metric space
 $({\mathcal A},\rho)$:
First, define $x\sim y$ if $\rho(x,y)=0$; clearly $\sim$ is an equivalence relation.
Denote the equivalent class containing $x$ by $[x]$. Set
$$\widetilde {\mathcal A}:=\mathcal A/\sim$$
to be the quotient space.
Secondly, for $[x],[y]\in\widetilde {\mathcal A}$, define
\begin{equation*}
 \widetilde{\rho}([x],[y])=\min\{\rho(a,b); a\in [x], b\in [y]\}.
\end{equation*}

\begin{thm}[\cite{HWYZ23}]\label{lem:induce}
The quotient space  $(\widetilde {\mathcal A}, \widetilde \rho)$ is a pseudo-metric space.
\end{thm}

Let $(\mathcal{A},\rho)$ be a pseudo-metric space. In the same manner as metric space, we can define
convergence of sequence, dense subset and completeness of ${\mathcal A}$.  (See \cite{Pepo90}.)
Recall that the H\"older and Lipschitz equivalence of pseudo-metric spaces are defined in \eqref{eq:Holder}.
The following theorem is obvious.

\begin{thm}[\cite{HWYZ23}]\label{extendLip}
Let $(\mathcal{A},\rho)$ and $(\mathcal{A}',\rho')$ be two complete pseudo-metric spaces.
Suppose  $B\subset\mathcal{A}$ is $\rho$-dense in $\mathcal{A}$ and  $B'\subset\mathcal{A}'$ is $\rho'$-dense in $\mathcal{A}'$. If $B\simeq B'$, then $\mathcal{A}\simeq\mathcal{A}'$.
\end{thm}

\medskip

\subsection{$\Sigma$-automaton}
~\\
\indent Recall that a \emph{finite state automaton (FSA)}  is  a 5-tuple
 $$(Q,\mathcal{A},\delta,q_0,P),$$
  where $Q$ is a finite set of states, $\mathcal{A}$ is a finite input alphabet, $q_0\in Q$ is the initial state, $P\subset Q$ is the set of final states, and $\delta$ is the transition function mapping $Q\times\mathcal{A}$ to $Q$. That is, $\delta(q,a)$ is a state for each state $q$ and input symbol $a$. (See for instance, \cite{JEH79}.)


Let $\Sigma=\{1,\dots, N\}$ be a finite set which we call an alphabet. For $a\in\Sigma$, we shall use $a^k$ to denote the word consisting of $k$ numbers of $a$.
Denote $\Sigma^{\infty}$ and $\Sigma^{k}$ to be the sets of infinite words and words of length $k$ over $\Sigma$, respectively.
Let $\Sigma^*=\bigcup_{k\geq0} \Sigma^{k}$ be the set of all finite words.

 \begin{defi}[\cite{HWYZ23}]\label{def:Sigma}
 \emph{A finite state automaton $M$ is called a \emph{$\Sigma$-automaton} if
\begin{equation*}
M=(Q,\Sigma^2,\delta,Id, Exit),
\end{equation*}
 where }

\emph{ (i) the state set is $Q=Q_0\cup \{Id, Exit\}$,  the initial state is $Id$, the final state is $Exit$;}

 \emph{(ii)  the input alphabet is $\Sigma^2$;}

 \emph{(iii)  the transition function $\delta$ satisfies
 \begin{equation*}
  \delta(Id,(i,j))=Id \Leftrightarrow i=j.
 \end{equation*}
 }
 \end{defi}

\medskip

\subsection{Surviving time}
~\\
\indent Now inputting symbol string $(\bx,\by)\in\Sigma^{\infty}\times\Sigma^{\infty}$ to $M$, we obtain a sequence of states $(S_{i})_{i\ge 0}$ and call it the \emph{itinerary} of $(\bx,\by)$. If we arrive at the state $Exit$, then we stop there and
 the itinerary is finite, otherwise,  it is infinite.
 We define the \emph{surviving time} of $(\bx,\by)$ to be
\begin{equation}\label{eq_surviviting_time}
T_M(\bx,\by)=\sup \{k;~ S_{k}\neq Exit\}.
\end{equation}

\begin{defi}
\emph{We say a $\Sigma$-automaton $M$ is \emph{feasible}, if there exists an integer $T_0\ge0$ such that
\begin{equation*}
  \min\{T_M(\bx,\by),T_M(\bx,\bz)\}\leq T_M(\by,\bz)+T_0, ~\forall \bx,\by,\bz\in\Sigma^{\infty}.
\end{equation*}
}
\end{defi}

\medskip

\subsection{Pseudo-metric space induced by $\Sigma$-automaton}
~\\
\indent Let $0<\xi<1$, we define a function $\rho_{M,\xi}$ on $\Sigma^\infty\times\Sigma^\infty$ as
\begin{equation*}
\rho_{M,\xi}(\bx, \by)=\xi^{T_{M}(\bx, \by)}.
\end{equation*}

If a $\Sigma$-automaton $M$ is feasible, then
$$\rho_{M,\xi}(\by, \bz)\leq \xi^{-T_0}\big(\rho_{M,\xi}(\bx, \by)+\rho_{M,\xi}(\bx, \bz)\big).$$
Hence $(\Sigma^{\infty}, \rho_{M,\xi})$ is a pseudo-quasi-metric space.
Let
\begin{equation*}
  (\cal A_M, \rho_{M,\xi})
\end{equation*}
 be the  pseudo-metric space  obtained from $(\Sigma^{\infty}, \rho_{M,\xi})$
 by the standard way in Section \ref{subsec:pseudo}, and we call it the \emph{pseudo-metric
space induced} by $M$.

\begin{lemma}[\cite{HWYZ23}]\label{Omegadense} Let $M$ be a feasible $\Sigma$-automaton and let  $(\mathcal{A}_M, \rho_{M,\xi})$
 be the pseudo-metric space induced by $M$. Let $\kappa\in \Sigma$. Then the set $\widetilde{\Omega}=\{[\omega\kappa^{\infty}];~\omega\in\Sigma^*\}$ is $\rho_{M,\xi}$-dense in $\mathcal{A}_M$.
\end{lemma}

\section{\textbf{Topology automaton of  Bara\'{n}ski carpet}} \label{Top-auto}

Let $K\in {\mathcal B}(n,m, \SD)$ be a Bara\'nski carpet generated by $\{\varphi_j\}_{j=1}^N$,  see \eqref{eq-IFS-B}. For $I=i_1\cdots i_k\in\Sigma^k$, denote
$\varphi_I=\varphi_{i_1} \circ \cdots\circ \varphi_{i_k}$ and we call $K_I=\varphi_I(K)$ a $k$-th cylinder.
We remark that in this section, we do not assume that $K$ satisfies the cross intersection condition.

\medskip

\subsection{Companion IFS}
~\\
\indent Let $K=K(n, m, \SD)$ be a Bara\'nski carpet, which  is the attractor of $\{\varphi_j\}_{j=1}^N$.
For each $\bd_j=(d_{j,1},d_{j,2})\in\SD$, define
$$
\tilde \varphi_j(x,y)= \left( \frac{x+d_{j,1}}{n}, \frac{y+d_{j,2}}{m} \right ).
$$
We call $\{\tilde \varphi_j\}_{j=1}^N$ the \emph{companion IFS} of $\{\varphi_j\}_{j=1}^N$ (or $K$).
Clearly the attractor of the companion IFS is either a fractal square, or a Bedford-McMullent carpet.

\medskip

\subsection{Topology automaton}
~\\
\indent Let $K\in {\mathcal B}(n,m,\SD)$ be a Bara\'nski carpet generated by the IFS $\{\varphi_j\}_{j=1}^N$.
Let $\{\tilde \varphi_j\}_{j=1}^N$ be the companion IFS of $\{\varphi_j\}_{j=1}^N$ and let
$K'$ be the corresponding attractor.
For $I,J\in\Sigma^k $, if $\tilde\varphi_I([0,1]^2)\cap\tilde\varphi_J([0,1]^2)\neq\emptyset$, then there are 8 possible positions between $\tilde\varphi_J([0,1]^2)$ and $\tilde\varphi_I([0,1]^2)$, which we will
indicate by elements in
 \begin{equation*}
Q_0=\{\pm \be_1, \pm \be_2, \pm(\be_1+\be_2), \pm(\be_1-\be_2)\}.
\end{equation*}

\begin{defi} [Topology automaton \cite{HWYZ23}]\label{def:topology}
{\rm
Let $\Sigma=\{1,\dots, N\}$.
We define   the \emph{topology automaton} of a fractal squre or a Bedford-McMullen carpet $K'$ to be the
 $\Sigma$-automaton
$$M_{K'}=\{Q_0\cup \{Id, Exit\}, \Sigma^2,  \delta, \{Id, Exit\}\}$$
  satisfying the following condition:
  for $i \neq j$ and $S\in Q_0\cup \{id\}$,
\begin{equation*}
\delta\big(S,(i,j)\big)=\left\{
\begin{array}{rl}
\tilde\varphi_i^{-1}(\tilde\varphi_{j}(\b0)+S),& \text{if~} \tilde\varphi_i(K)\cap (\tilde\varphi_{j}(K)+\bs)\neq \emptyset,\\
Exit,& \text{otherwise},
\end{array}
\right.
\end{equation*}
where we regard $id$ as $\b0$ in the right hand side of the above formula.
}
\end{defi}


\begin{defi}\emph{We define the topology automaton of a Bara\'nski carpet $K$
to be the automaton of $K'$, that is,
$M_{K}:=M_{K'}.$
}
\end{defi}

For $\bx=x_1x_2\dots\in\Sigma^\infty$, denote by $\bx|_k=x_1\dots x_k$ the prefix of $\bx$ with length $k$.

\begin{thm}\label{thm:NA}
 The topology automaton  $M_K$ of a Bara\'nski carpet $K$  is a feasible $\Sigma$-automaton.
 Precisely,
 \begin{equation}\label{eq:xyz}
 \min\{T(\bx,\by), T(\bx,\bz)\}\leq T(\by,\bz)+1,\quad \forall \bx,\by,\bz\in \Sigma^\infty.
 \end{equation}
\end{thm}

\begin{proof} Without loss of generality, we assume that $K$ is a fractal square or a Bedford-McMullen carpet.
Let $\bx, \by, \bz$ be three different points in $\Sigma^\infty$.
Denote $k=T(\by,\bz)$. Then $k$ is the smallest integer such that
$K_{\by|_k}\cap K_{\bz|_k}=\emptyset$.

Suppose on the contrary that \eqref{eq:xyz} does not hold. Then
   $\varphi_{\bx|_{k+2}}([0,1]^2)$ intersects both  $\varphi_{\by|_{k+2}}([0,1]^2)$ and $\varphi_{\bz|_{k+2}}([0,1]^2)$,
   which implies that  $\varphi_{\by|_{k}}([0,1]^2)\cap \varphi_{\bz|_{k}}([0,1]^2)\neq \emptyset$.

\textit{Case 1.}   $\varphi_{\by|_k}(K)$ and $\varphi_{\bz|_k}(K)$ are located in the same row (or column) of oder $k$.

Since $\varphi_{\bx|_k}([0,1]^2)$ intersects both  $\varphi_{\by|_k}([0,1]^2)$ and $\varphi_{\bz|_k}([0,1]^2)$,
without loss of generality, we may assume that $\varphi_{\bx|_k}([0,1]^2)$  locates under $\varphi_{\by|_k}([0,1]^2)$ and $\varphi_{\bz|_k}([0,1]^2)$ locates on the right side of $\varphi_{\by|_k}([0,1]^2)$, see Figure \ref{fig:XYZ}.
This forces
$$x_{k+1}=(n-1, m-1)\in \SD \text{ \ and  \ } z_{k+1}=(0,0)\in \SD.$$
(See the blue rectangle and red rectangle in Figure \ref{fig:XYZ}.)
We argue that $y_{k+1}\neq (n-1,0)$, for otherwise,  $\varphi_{\bz|_{k}}(0)\in \varphi_{\by|_k}(K)\cap \varphi_{\bz|_k}(K)$,  which contradicts the maximality of $k$.
(See the yellow rectangle in Figure \ref{fig:XYZ}.)

\begin{figure}[H]
	\centering
	\includegraphics[width=7.5 cm]{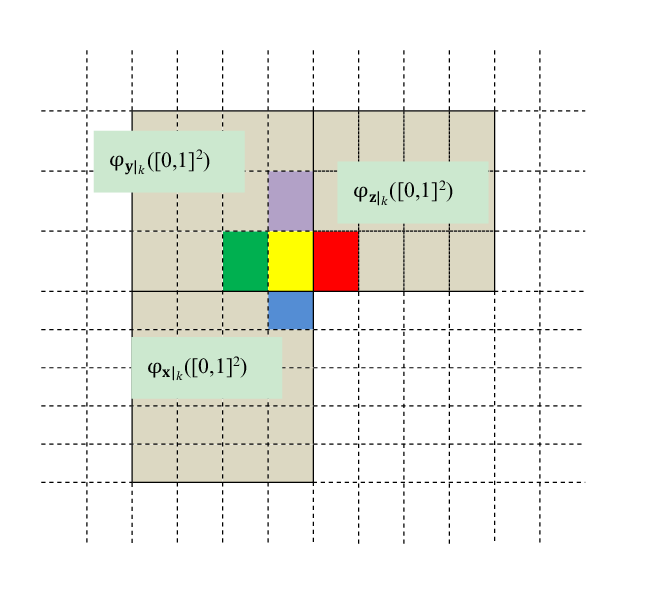}\\
	\caption{Illustration of proof of Theorem \ref{thm:NA}: Case 1.}\label{fig:XYZ}
\end{figure}

To guarantee
$\varphi_{\bx|_{k+2}}[0,1]^2\cap   \varphi_{\bz|_{k+2}}[0,1]^2\neq \emptyset$, we must have
$x_{k+2}=(n-1, m-1)$, but then
  $\varphi_{\bx|_{k+2}}([0,1]^2)$  and
  $\varphi_{\by|_{k+2}}([0,1]^2)$ are disjoint,  which contradicts to our assumption.

 \smallskip
\textit{Case 2.}  $\varphi_{\by|_k}([0,1]^2)\cap \varphi_{\bz|_k}([0,1]^2)$ is a single point.

By a similar argument as Case 1, one can show that \eqref{eq:xyz} holds in this case.
\end{proof}

\medskip

\subsection{Coding and projection}
~\\
\indent Let $K$ be  the attractor of an IFS $\{\varphi_j\}_{j=1}^N$.
Define $\pi_{K}:\Sigma^{\infty}\to K$, which we call a \emph{projection}, by
$$
\big\{\pi_{K}(\bx)\big\}=\bigcap_{k\geq1} \varphi_{x_1\cdots x_k}(K).
$$
If $\pi_{K}(\bx)=x\in K$, then we call the sequence $\bx=x_1x_2\dots\in\Sigma^\infty$ a \emph{coding} of $x$.

Let $\xi\in (0,1)$.
Let $(\cal A_{M_K}, \rho_{M_K,\xi})$ be the pseudo-metric space induced by $M_K$, see Section 2.4.
The following lemma is obvious.

\begin{lemma}
That  $\pi_K(\bx)=\pi_K(\by)$ if and only if $\rho_{M,\xi}(\bx,\by)=0$.
\end{lemma}

We define $\pi: \cal A_{M_K}\to K$ by
$$
\pi([\bx])=\pi_K(\bx).
$$

For $\varphi_j=(f_{d_{j,1}},g_{d_{j,2}})$, we denote $\varphi'_j=(f'_{d_{j,1}},g'_{d_{j,2}})=(a_j, b_j)$, and denote
$$
r^*=\max\{a_j, b_j; ~1\leq j\leq N\}, \quad r_*=\min\{a_j, b_j; ~1\leq j\leq N\}.
$$

\begin{lemma}\label{lem:C3} There is a constant $C_3>0$ such that
if $I,J\in \Sigma^k$ and $K_I\cap K_J=\emptyset$, then
$\text{dist}(K_I, K_J)>C_3r_*^k$.
\end{lemma}

\begin{proof} If $\varphi_I([0,1]^2\cap \varphi_J([0,1]^2)=\emptyset$, then clearly
$\text{dist}(K_I, K_J)\geq r_*^k$. So in the following, we assume that $\varphi_I([0,1]^2\cap \varphi_J([0,1]^2)\neq \emptyset$.
Let
$$\delta_0=\min\left ( \{\text{dist~}(K, K+\bb);~\bb\in \pm\{\be_1, \be_2, \be_1+\be_2, \be_1-\be_2\}, K\cap (K+\bb)\neq \emptyset\}\right ).$$
We will show that
\begin{equation}\label{dist-KK}
\text{dist}(K_I, K_J)\geq \delta_0 r_*^k.
\end{equation}

Recall that $\varphi_I'=(a_I, b_I)$, $\varphi_J'=(a_J,b_J)$.

\textit{Case 1.}   $\varphi_I([0,1]^2)$ and    $\varphi_J([0,1]^2)$ locate in the same row or in the same column.

Without loss of generality, let us assume that $K_I$ and $K_J$ locate in the same row, then $b_I=b_J$.
Again, without loss of generality, let us assume  $a_I\leq a_J$.
Notice that $K_I\cap K_J=\emptyset$ implies that $K\cap (K+\be_1)=\emptyset$.
Since
$$\varphi_I^{-1}\varphi_J(K)=\be_1+\text{diag}\left (\frac{a_J}{a_I},1\right ) \cdot K, $$
 we obtain $\text{dist}(K, \varphi_I^{-1}\varphi_J(K))\geq \text{dist}(K,  K+\be_1)\geq\delta_0$.
Consequently, \eqref{dist-KK} holds.

\textit{Case 2.} $\varphi_I([0,1]^2)$ and    $\varphi_J([0,1]^2)$ meet at a corner.

Without loss of generality, let us assume that the most right-bottom point of $\varphi_I([0,1]^2)$ coincides with the most left-top point of $\varphi_J([0,1]^2)$.
Denote the intersection point by $z$.
Let $f(x)=r_*^{-k}(x-z)$.
Then $f(K_I)$ and $f(K_J)$ meet at $0$, and
$$f(K_I)=\text{diag}\left (\frac{a_I}{r_*^k}, \frac{b_I}{r_*^k}\right )\cdot (K-\be_1), \quad
f(K_J)=\text{diag}\left ( \frac{a_J}{r_*^k}, \frac{b_J}{r_*^k}\right )\cdot (K-\be_2),$$
Therefore, $\text{dist~}(f(K_I),f(K_J))\geq \text{dist}(K, K+(\be_1-\be_2))\geq \delta_0$, and  \eqref{dist-KK} follows.
The lemma is proved.
\end{proof}

\begin{theorem}\label{thm: Holder map}
Let $K$ be a  Bara\'{n}ski carpet,
Let $s=\sqrt{\log r^*/\log r_*}$ and $\xi=(r_*)^s$.
Then $\pi: (\cal A_{M_K}, \rho_{M_K,\xi})\to K $ is a bi-H\"{o}lder map with index $s$.
In particular, if $K$ is a fractal square, then $\pi$ is bi-Lipschitz.
\end{theorem}

\begin{proof}
Take $x,y\in K$. Let $\bx$  be a coding of $x$ and $\by$ be a coding of $y$.
Let $k = T(\bx, \by)$ be the surviving time of $(\bx,\by)$ in the  automaton $M_K$, see \eqref{eq_surviviting_time}.
Then the $k$-th cylinder containing $x$ and that containing $y$, either coincide or have non-empty intersection.
Since every $k$-th cylinder has diameter no larger than $2(r^*)^k$, it follows  that
$$
|x-y|\leq 4(r^*)^k.
$$

On the other hand, the $(k+1)$-th cylinder containing $x$ and that containing $y$ are disjoint,
so   we have
$$
|x-y|\geq C_3(r_*)^{k+1},
$$
where $C_3$ is the constant in Lemma \ref{lem:C3}.
Notice that
\begin{equation*}
  \rho_{M_K,\xi}(\bx,\by)=\xi^k=(r_*)^{sk}=(r^*)^{k/ s}.
\end{equation*}
Set $C= \max\{4, 1/(C_3r_*)\}$, we obtain the theorem.
\end{proof}

\medskip

\subsection{Symmetric $\Sigma$-automaton}
~\\
\indent Let $M=(\Sigma, Q, \delta, Id, Exit)$ be a $\Sigma$-automaton.
We say $M$ is \emph{symmetric}, if

(i) $Q_0=Q_1\cup \overline{Q}_1$ with $Q_1\cap \overline{Q}_1=\emptyset$,
and there exist a bijection $\varrho:Q_1\to \overline{Q}_1$.
We call $\varrho(S)$  the \emph{mirror state} of $S$.
For simplicity, hereafter we denote $\varrho(S)$ by $-S$.

(ii) By convention, we set $-Id=Id$ and $-Exit=Exit$.

(iii) For any $\bs\in Q$, $i,j\in\Sigma$, it holds that
$$
\delta\big(\bs,(j,i)\big)=-\delta\big(-\bs,(j,i)\big).
$$

Clearly, the topology automaton of a Bara\'nski carpet is symmetric.

\section{\textbf{Cross automaton}}\label{CA}

To study  Bara\'nski carpets satisfying the cross intersection condition,
we introduce  the cross automaton as following.

\begin{defi}[Cross automaton]\label{def:cross auto}
{\rm A symmetric $\Sigma$-automaton $M=\{Q,\Sigma^2,\delta,Id,Exit\}$ is called a \emph{cross automaton} if
 $Q_0=\{\pm \be_1, \pm \be_2\}$
and  $M$ satisfies the following  conditions:

(i) (Uniqueness) If $\delta\big(S,(i,j_1)\big)=\delta\big(S,(i,j_2)\big)\neq Exit$, then $j_1=j_2$.

(ii) (Self-looping property) For any $S\in Q_0$ and $(i,j)\in\Sigma^2$,
$\delta\big(S,(i,j)\big)\in\{S,Exit\}$. (That is,
 a state $S$ either transfers to itself  or  to $Exit$.)

(iii) (Triple-coding-free condition) Let $\bx,\by,\bz\in\Sigma^{\infty}$ be distinct, then
at most one of $T_M(\bx,\by),T_M(\bx,\bz)$ and $T_M(\by,\bz)$ take the value $+\infty$.
}
\end{defi}

\begin{figure}[h]
  \centering
  \includegraphics[width=10cm]{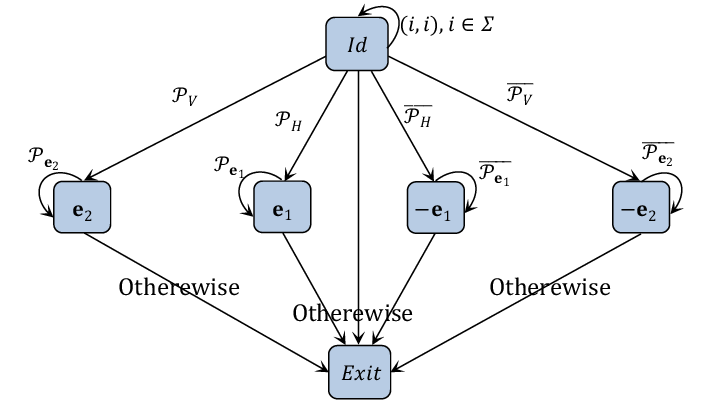}
  \caption{The transition diagram of a   cross automaton.}
  \label{diagram}
\end{figure}

For $S,S'\in Q_0\cup\{Id\}$, we define
$$
\mathcal{P}_{S\to S^{'}}=\left\{(i,j)\in\Sigma^2;\delta\big(S,(i,j)\big)=S^{'}\right\}.
$$
For convenience,  we denote
$$
  \cal P_H=\mathcal{P}_{Id\to \be_1},\quad \cal P_V=\mathcal{P}_{Id\to \be_2},\quad
  \cal P_{\be_1}=\mathcal{P}_{\be_1\to \be_1},\quad \cal P_{\be_2}=\mathcal{P}_{\be_2\to \be_2}.
$$
For a set $\cal P\subset\Sigma^2$, we define $\overline{\cal P}=\{(j,i);(i,j)\in \cal P\}$.

Figure \ref{diagram} illustrates the transition diagram of a  cross automaton $M$.
Clearly, $M$ is completely determined by the sets $\mathcal{P}_{H}$, $\mathcal{P}_{V}$, $\mathcal{P}_{\be_1}$ and $\mathcal{P}_{\be_2}$.

If $(i, j)\in\cal P_H$, then we denote $i\lhd_H j$ and say $i$ is the \emph{H-predecessor} of $j$ and $j$
is the \emph{H-successor} of $i$; similarly, we can define $i\lhd_V j$, $i\lhd_{\be_1} j$ and $i\lhd_{\be_2} j$.

The following lemma is obvious.

\begin{thm} Let $K$ be a Bara\'nski carpet satisfying the cross intersection condition, then
the topology automaton $M_K$ is a cross automaton. Moreover, if $K$ satisfies the vertical separation condition,
then   ${\cal P}_V=\emptyset$.
\end{thm}

\begin{proof} We only need to show the theorem holds for a fractal square or a Bedford-McMullen carpet.
By the cross intersection condition, the states $\pm (\be_1+\be_2)$ and $\pm(\be_1-\be_2)$ will not occur
in $M_K$, so $Q_0=\{\pm \be_1, \pm \be_2\}$.

Now we verify item (i)-(iii) in Definition \ref{def:cross auto}.
Item (i) holds since for every $k$-th cylinder $K_I$, there is at most one $k$-th cylinder
locates on the right (or left, or above, or below) of $K_I$ and adjacent to $K_I$.
Item (ii) holds since if two $k$-th cylinders $K_I$ and $K_J$ locate in the same row (\textit{resp.}  column), then
$K_{Ii}$ and $K_{Jj}$ have no chance to locate in the same column (\textit{resp.} row).
Item (iii) holds since by the cross intersection condition, no points of $K$ has more than two codings.
\end{proof}

The following example shows that the class of cross automata is much wider than that of topology automata of Bara\'nski carpets.

\begin{figure}[H]
  \centering
  \subfigure[]{
   \includegraphics[width=5.2 cm]{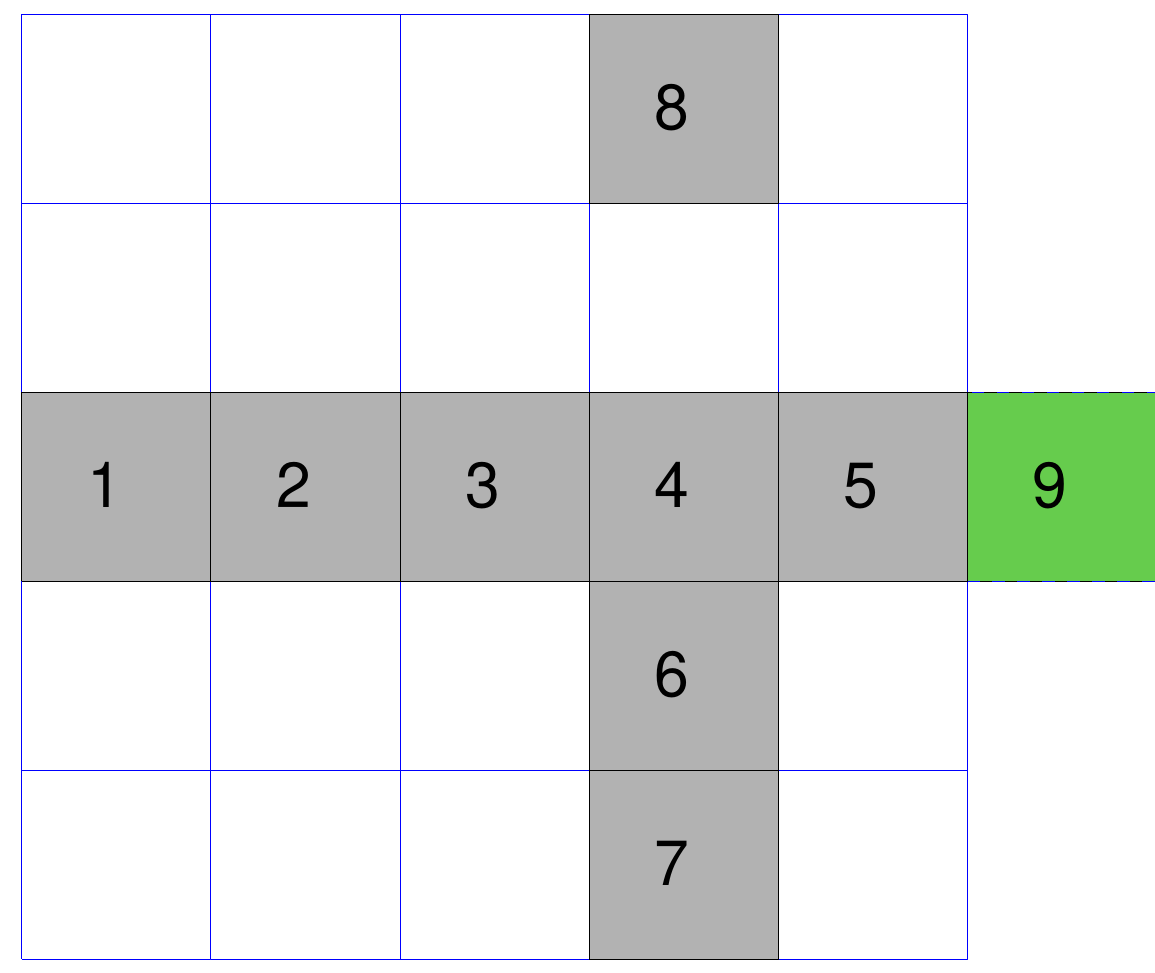} }
   \subfigure[]{
  \includegraphics[width=6.5 cm]{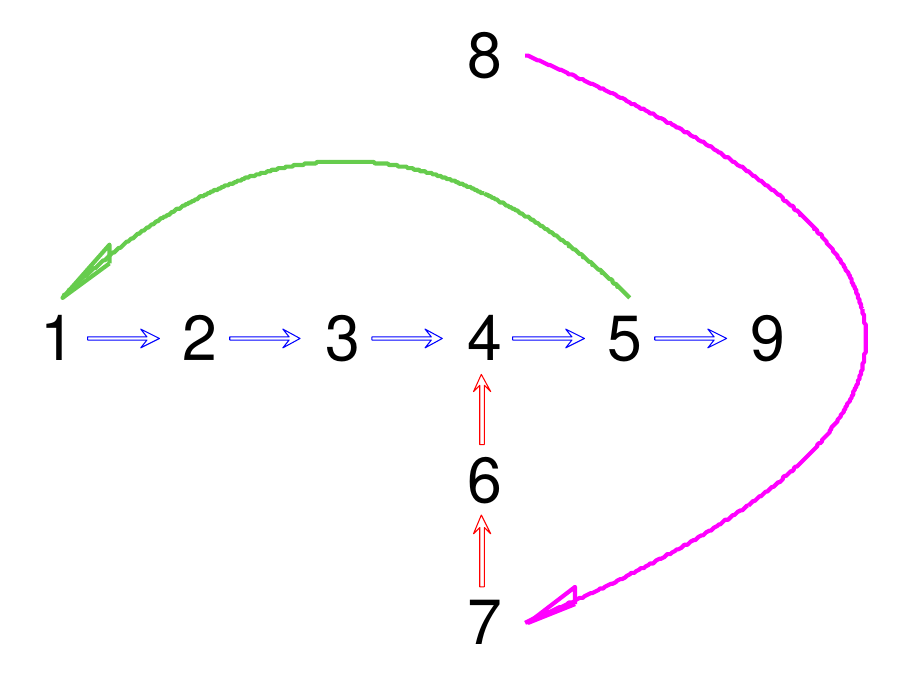}
  }

  \caption{(a) The cross automaton in Example 4.1.  (b) The graph representation of the cross automaton in (a).
  (We draw the four graphs in one picture, but we use edges of different colors to distinct them.)}
  \label{fig:graph}
\end{figure}

\begin{example} \emph{ Let $K$ be the fractal square consisting of $8$ cylinders
indicated in Figure \ref{fig:graph}(a). Let $M_K$ be the topology automaton of $K$.
Let $M$ be a $\Sigma$-automaton with    $\Sigma=\{1,\dots, 9\}$
 such that $M$ is an extension of $M_K$, that is, all edges of $M_K$ belong to $M$.
Moreover, we set
$$\delta(Id, (9,9))=Id,  \ \  \delta(Id, (5, 9))=\be_1, \ \ \delta(Id, (9,5))=-\be_1,$$
and set $\delta(S, (i,9))=\delta(S, (9,i))=Exit$ otherwise.  Precisely, we have
$$
{\cal P}_H=\{(1,2),(2,3),(3,4),(4,5),(5,9)\},
$$
$$
{\cal P}_V=\{(7,6),(6,4)\}, \quad {\cal P}_{\be_1}=\{(5,1)\},\quad {\cal P}_{\be_2}=\{(8,7)\}.
$$
Then $M$ is a cross automaton but it
is not the topology automaton of any Bara\'nski carpet.
}
\end{example}


\medskip

\subsection{Graph representation of  cross automaton}
~\\
\indent Firstly, we recall some notions of graph theory, see \cite{Bal2000}. Let $G=(V, \mathcal{E})$ be a directed graph, where $V$ is the \emph{vertex set} and $\mathcal{E}$ is the \emph{edge set}. Each edge $\mathbf{e}\in \cal E$ is associated to an ordered pair $(u,v)$ in $V\times V$, and we say
$\mathbf{e}$ is \emph{incident out} of $u$ and \emph{incident into} $v$. Denote that $\be=(u,v)$.
We also call $u$ and $v$ the \emph{origin} and \emph{terminus} of $\mathbf{e}$, respectively. The number of edges incident out of a vertex $v$ is the \emph{outdegree} of $v$ and is denoted by $\deg^-(v)$. The number of edges incident into a vertex $v$ is the \emph{indegree} of $v$ and is denoted by $\deg^+(v)$.
If $\deg^-(v)=0$, then we say $v$ is \emph{maximal}, if $\deg^+(v)=0$, then we say $v$ is \emph{minimal}. If $v$ is both minimal and maximal, then we say $v$ is \emph{isolated}.

A \emph{directed walk} joining vertex $v_1$ to vertex $v_k$ in $G$ is a sequence $(v_1, v_2, \dots, v_k)$ with $(v_{i},v_{i+1})\in\mathcal{E}$, in addition, if all $v_i (1\le i\le k)$ are distinct, then we call it a \emph{path}. If all $v_i (1\le i\le k-1)$ are distinct and $v_k=v_1$, then we call it a \emph{cycle}. Let $(v_1, v_2, \dots, v_k)$ be a path, if $v_1$ is minimal and $v_k$ is maximal, then we call it a \emph{chain}.

For a  cross automaton $M$, we will regard $(\Sigma,\cal P_H)$   a graph. Precisely, there is an edge in $(\Sigma,\cal P_H)$ from $i$ to $j$
if and only if $(i,j)\in {\mathcal P}_H$.
A symbol $j\in\Sigma$ is said to be \emph{H-minimal} (resp. \emph{maximal}) if it is minimal (resp. maximal) in $(\Sigma,\cal P_H)$.

Similarly, we   define $(\Sigma,\cal P_V)$, $(\Sigma,\cal P_{\be_1})$, $(\Sigma,\cal P_{\be_2})$
as well as \emph{V-minimal} (\emph{maximal}),  \emph{ $\be_1$-minimal} (\emph{maximal})
and \emph{ $\be_2$-minimal} (\emph{maximal}). One can refer Figure \ref{fig:graph} as an example.

\subsection{Feasibility of cross automaton}
~\\
\indent We write the initial state by $id$ instead of $Id$ for clarity.  we will use $S\stackrel{(i,j)}\longrightarrow S'$ as an alternative notation for $\delta(S,(i,j))=S'$.

The following theorem  says that  all cross automata are  feasible.

\begin{thm}\label{dis}Let $M$ be a  cross automaton.
	For $\bx,\by,\bz\in\Sigma^\infty$, we have
	\begin{equation}\label{dis1}
		\min\{T_M(\bx,\by), T_M(\bx,\bz)\}\le T_M(\by,\bz)+1.
	\end{equation}
\end{thm}
\begin{proof}
Clearly, \eqref{dis1} holds if any two of $\bx,\by,\bz$ are identical, or $T_M(\by,\bz)=\infty$.
Now we assume $\bx,\by,\bz$ are distinct and
	  $q=T_M(\by,\bz)<\infty$.

 Suppose on the contrary that \eqref{dis1} is false.   Then $T_M(\bx,\by)\ge q+2$ and $T_M(\bx,\bz)\ge q+2$.

Denote the itinerary of $(\bx,\by)$, $(\bx,\bz)$ and $(\by,\bz)$  by $(S_k)_{k\geq 0}$,
  $(S'_k)_{k\geq 0}$ and  $(U_k)_{k\geq q+1}$, respectively.
  Here we have   $S_0=S'_0= U_0=id$. Notice that   $U_{q+1}=Exit$.

\textit{Case 1.}  $S_{q+1}=Id$ or $S'_{q+1}=Id$.

Without loss of generality, let us assume that $S_{q+1}=Id$.
This means that $\bx|_{q+1}=\by|_{q+1}$. So
 $(S'_k)_{k=0}^{q+1}=(U_k)_{k=0}^{q+1}.$
Since $S'_{q+1}\neq Exit$, we deduce that $U_{q+1}\neq Exit$, which contradicts to $T_M(\by,\bz)=q$.

\textit{Case 2.}  $S_{q+1}\neq Id$ and $S'_{q+1}\neq Id$.

Dentoe  $S_{q+1}=\bs$ and $S'_{q+1}=\bs'$. By the self-looping property, we have $S_{q+2}=\bs$, so
 $$\bs \stackrel{(x_{q+2},y_{q+2})}\longrightarrow \bs.$$
 By the same reason, we have
 $\bs'\stackrel{(x_{q+2},z_{q+2})}\longrightarrow \bs'.$  Set
$$
\bx'=(x_1\dots x_{q+1})(x_{q+2})^\infty, \quad  \by'=(y_1\dots y_{q+1})(y_{q+2})^\infty, \quad
\bz'=(z_1\dots z_{q+1})(z_{q+2})^\infty,.
$$
Then the itinerary of $(\bx',\by')$ is $id\mapsto (S_1\dots S_{q})(\bs)^\infty$, and it follows $T(\bx',\by')=\infty$. Similarly, $T(\bx',\bz')=\infty$.
This   contradicts  to the  triple-coding-free  condition. 

The contradictions in the two cases prove the theorem.   
\end{proof}

\section{\textbf{Simplification of cross automaton}}\label{one Sim}

In this section, we study the simplification of cross automaton.

\subsection{Cross automaton of Class 0, Class 1, and Class 2}
~\\
\indent We will confine ourself to three special classes of cross automata.
\begin{defi}
\emph{
Let $M$ be a cross automaton.}
	
	\emph{(i) We say $M$ is  \emph{of Class $0$} if   ${\mathcal P}_V=\emptyset$. }
	
	\emph{(ii) We say $M$ is \emph{of Class $1$} if $M$ is the topology automaton
		of a Bara\'nski carpet satisfying  the cross intersection condition as well as the top isolated condition, but does not satisfy
		the vertical separation condition.}
\end{defi}

 Let $K$ be a Bara\'nski carpet satisfying the cross intersection condition as well as the top isolated condition,
but does not satisfy the vertical separation condition.
Let $\gamma\in \Sigma$ be the letter such that $\varphi_\gamma(K)$ is the top cylinder of $K$.
There exists $\lambda\in \Sigma$ such that $\varphi_\lambda(K)$ locates in the same column as $\varphi_\gamma(K)$ and in the bottom of $[0,1]^2$.

\begin{lemma} \label{def:gamma-isolated}
	Let $M$ be a   cross  automaton of Class $1$. Then

\indent\emph{(i)} There exists $\gamma\neq \lambda\in \Sigma$ such that $\cal P_{\be_2}=\{(\gamma,\lambda)\}$; we call $\gamma$ the \emph{top vertex}  and $\lambda$ the \emph{bottom vertex} of $M$. \\	
\indent\emph{(ii)} $\gamma$ is $H$-isolated, $V$-isolated and  $\be_1$-isolated.\\
	\indent\emph{(iii)}   For any $\theta_1\neq \lambda$ and $\theta_2\in \Sigma$, the inputs
$(\lambda\lambda, \theta_1\theta_2)$ and $(\theta_1\theta_2,\lambda\lambda)$ lead $id$ to $Exit$. \\
	\indent\emph{(iv)} ${\mathcal P}_V\neq \emptyset$ and the graph $(\Sigma,\cal P_V)$ has no cycle.
\end{lemma}
\begin{proof} (i) and (ii) are obvious.
Item (iii)  means that $K_{\lambda\lambda}$ does not intersect other
cylinders $K_\theta$. If $K_\lambda$ intersects the right boundary of $[0,1]^2$, then $K_{\lambda\lambda}$
locates on the most right-bottom corner of $[0,1]^2$, so it does not intersect any other cylinder $K_\theta$.
If $K_\lambda$ intersects the left boundary of $[0,1]^2$, the result holds by the same argument.
If $K_\lambda$ does not intersect the right and the left boundaries of $[0,1]^2$, then  $K_{\lambda\lambda}$ does not intersect the right and  the left boundaries of $K_\lambda$, so the result is also true.

	
	(iv) That $(a,b)\in {\mathcal P}_V$ implies that $\phi_b ([0,1]^2)$ is adjacent and above $\phi_a ([0,1]^2)$, so there
	is no cycle.
\end{proof}

\begin{defi} \emph{Let $M$ be a cross automaton and let $\{\gamma, \lambda\}\in \Sigma$.
		If all the items (i)-(iv) in Lemma \ref{def:gamma-isolated} hold for $M$, then we call $M$ a cross automaton of \emph{Class 2}. }
\end{defi}

Clearly, Class 1 is a sub-family of Class 2, and the later one  is much wider.

\medskip
\subsection{Simplification of  cross automaton}
~\\
\indent Let $M$ be a cross automaton of Class $2$.
  Since $\cal P_V\ne\emptyset$ and the graph $(\Sigma,\cal P_V)$ has no cycle,  we can pick  $(\tau,\kappa)\in\mathcal{P}_{V}$ such that $\kappa$ is $V$-maximal.

Let
\begin{equation}\label{fig_break1}
	\mathcal{P}'_{V}=\mathcal{P}_{V}\setminus\{(\tau,\kappa)\},~
	\mathcal{P}'_{H}=\mathcal{P}_{H},~
	\mathcal{P}'_{\be_1}=\mathcal{P}_{\be_1}~{\rm and}~
	\mathcal{P}'_{\be_2}=\mathcal{P}_{\be_2}.
\end{equation}
Let $M'$ be the  cross automaton determined by $\mathcal{P}'_{H},\mathcal{P}'_V$, $\mathcal{P}'_{\be_1}$ and $\mathcal{P}'_{\be_1}$,
and we call it the \emph{one-step simplification} of $M$ by deleting $(\tau,\kappa)$.

\begin{lemma}\label{M'carpetauto}
 Let $M$ be a   cross automaton of Class 2.
	If $M'$ is a one-step simplification of $M$ by deleting $(\tau,\kappa)$ in the graph $(\Sigma,\mathcal{P}_{V})$. Then we have
	
	\emph{(i)} $M'$ is a  cross automaton of   Class 0  if $\mathcal{P}'_V=\emptyset$, and is of  Class 2 otherwise;
	
	\emph{(ii)} $\tau$ is $V$-maximal and $\kappa$ is $V$-isolated in   $M'$.
\end{lemma}
\begin{proof}
	(i) Notice that $(\Sigma,\mathcal{P}'_{H})=(\Sigma,\mathcal{P}_{H})$,
	$(\Sigma,\mathcal{P}'_{\be_i})=(\Sigma,\mathcal{P}_{\be_i}) (i=1,2)$,
	and $(\Sigma,\mathcal{P}'_{V})$ is a subgraph of $(\Sigma,\mathcal{P}_{V})$,
	we infer that $M'$ satisfies uniqueness property,
	the self-looping property and triple-coding-free condition in the definition of cross automaton.
    By the same reason,   $M'$  satisfies   item (i)-(iv) in Lemma \ref{def:gamma-isolated}.
	
	(ii) That  $\kappa$ is $V$-maximal in $M$ means $\kappa$ is maximal in $(\Sigma,\mathcal{P}_{V})$ and thus also is maximal in $(\Sigma,\mathcal{P}'_{V})$, so $\kappa$  is  $V$-maximal in $M'$.
	The process of one-step simplification breaks one edge of $(\Sigma,\mathcal{P}_{V})$  from $\tau$ to $\kappa$,
	therefore, by uniqueness property, $\tau$ is $V$-maximal in $M'$ and $\kappa$ is $V$-isolated in $M'$.
\end{proof}

\begin{figure}[H]
	\centering
	\subfigure[$K_1$: before simplification.]{
		\includegraphics[width=4.5cm]{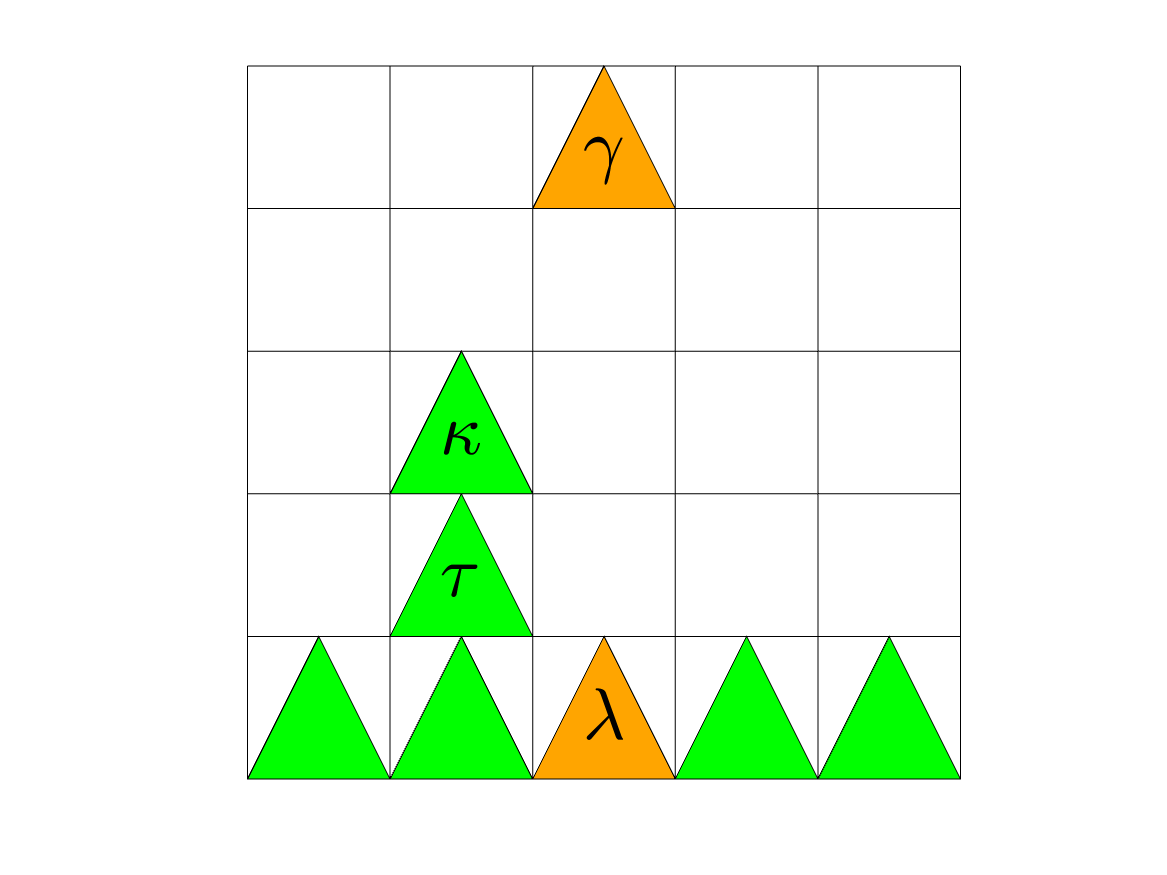}
	}
	\subfigure[$K_2$: after simplification.]{
		\includegraphics[width=4.5cm]{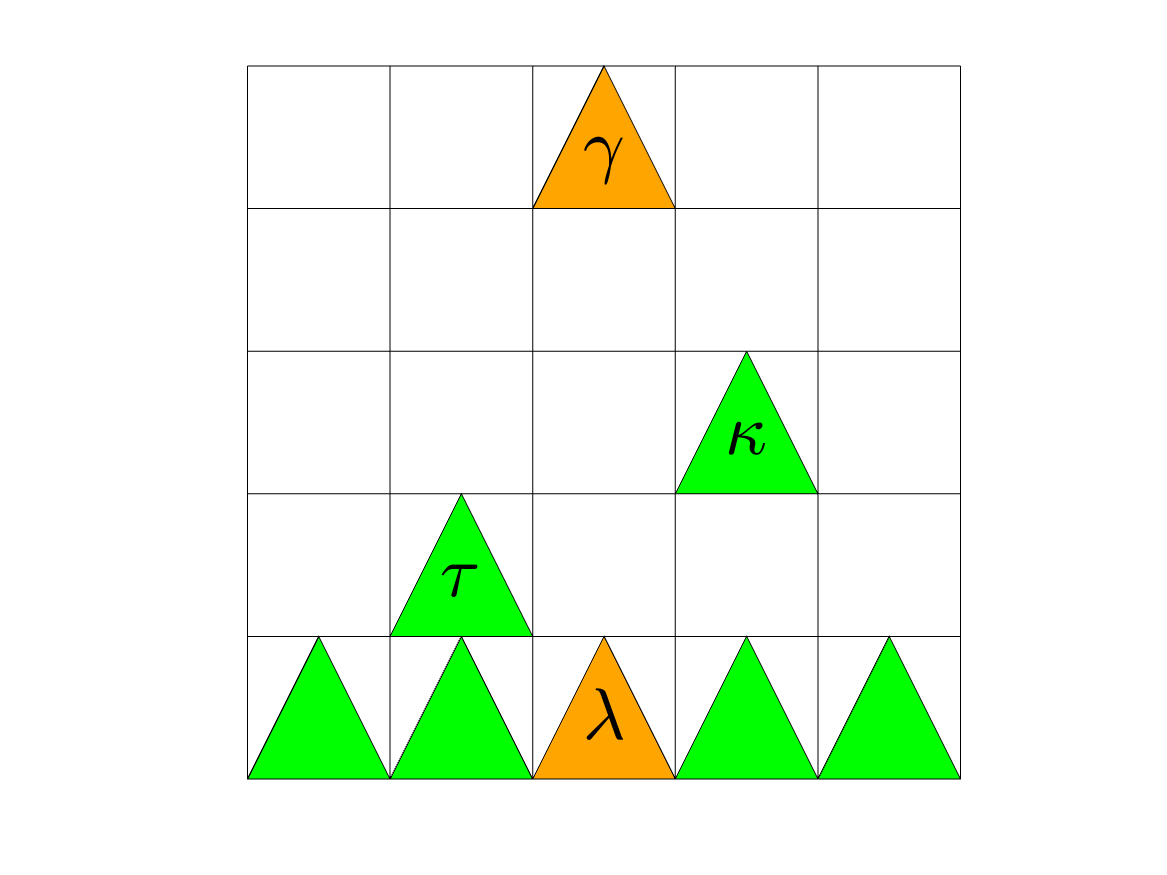}
	}\\
	\subfigure[$K_2$: before simplification.]{
		\includegraphics[width=4.5cm]{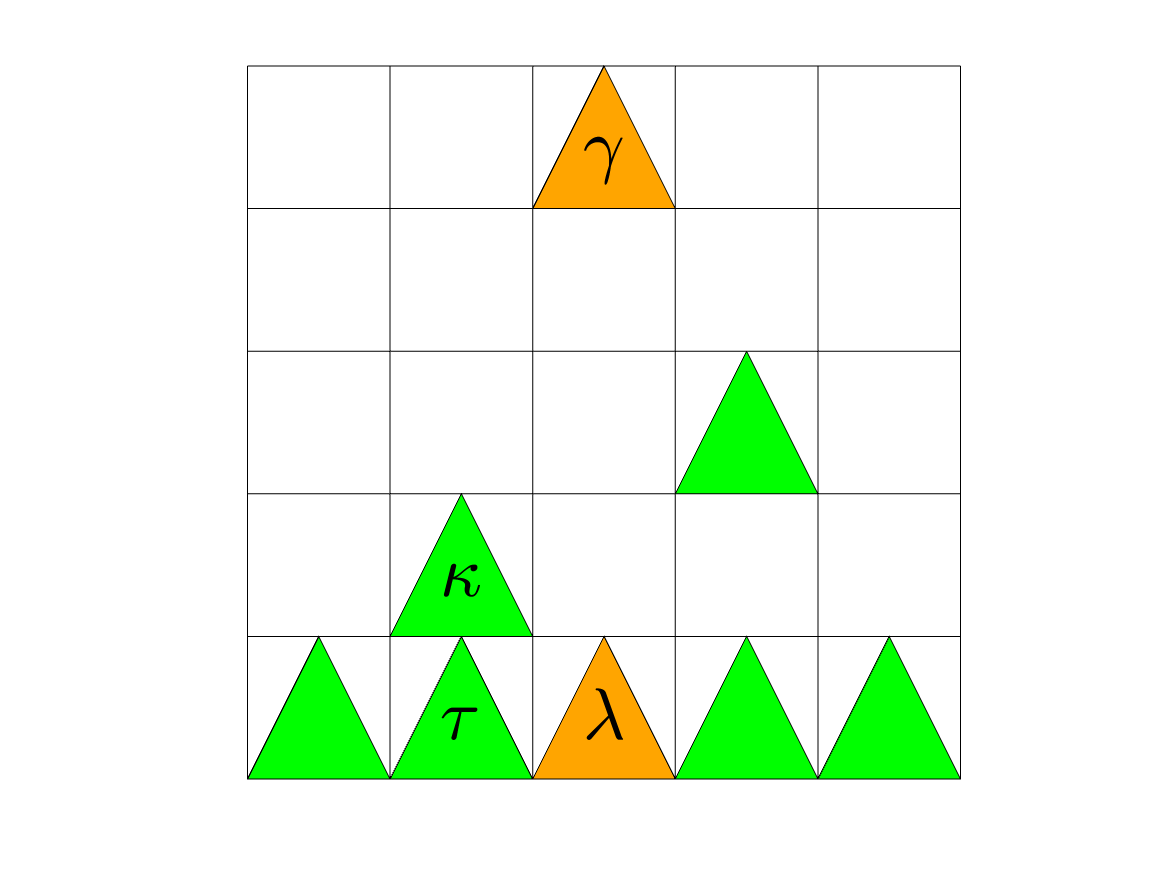}
	}
	\subfigure[$K_3$: after simplification.]{
		\includegraphics[width=4.5cm]{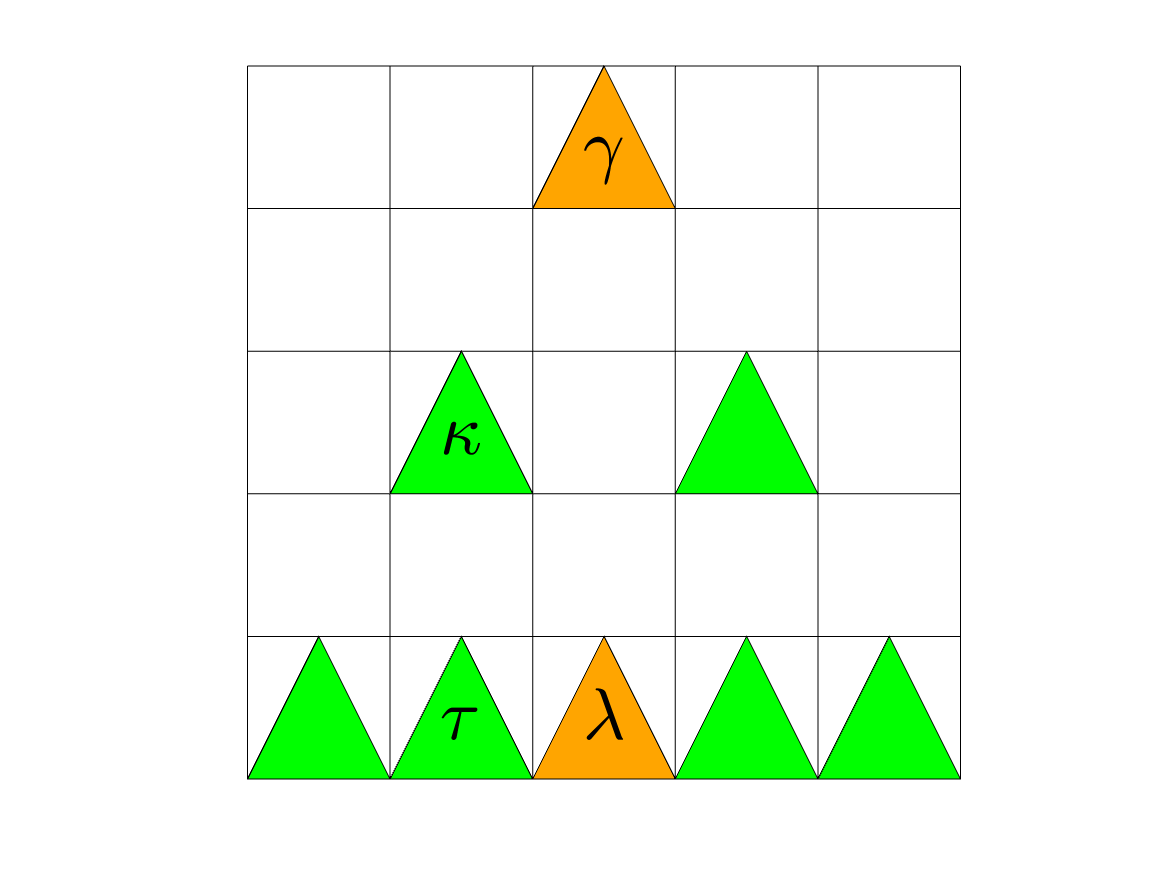}
	}\\
	\caption {$M_{K_2}$ is a one step simplification of $M_{K_1}$, and $M_{K_3}$ is a one step simplification of $M_{K_2}$.}
	\label{fig:6}
\end{figure}

\begin{thm}\label{main} Let $M$ be a cross automaton of Class 2, and let $M'$ be a one-step simplification of $M$. Then for any $\bx,\by\in\Omega=\{\omega \kappa^{\infty};\omega\in\Sigma^*\}$ there exists a bijection $g:\Omega\to \Omega$ such that
	\begin{equation}\label{eq-key}
		|T_M(\bx,\by)-T_{M'}(g(\bx),g(\by))|\leq 4.
	\end{equation}
\end{thm}

Theorem \ref{main} is the most important theorem in this paper, and the proof is rather technical.
We will prove it in Section \ref{pf main}. As a corollary of Theorem \ref{main}, we have

\begin{thm}\label{thm:onesimp}
	Let $M$ be a  cross automaton of Class 2 and
	let $M'$ be a one-step simplification of $M$. Then $(\cal A_M,\rho_{M,\xi})\simeq (\cal A_{M'},\rho_{M',\xi})$ for any $\xi\in(0, 1)$.
\end{thm}
\begin{proof}
	Recall that $\Omega=\{\omega\kappa^{\infty};\omega\in\Sigma^*\}.$ Let $\widetilde{\Omega}=\{[\omega\kappa^{\infty}];\omega\in\Sigma^*\}$.
	Let $g:\Omega\rightarrow\Omega$ be the bijection map given in Theorem \ref{main}. Since for any $\bx,\by\in\Omega$,
	\begin{equation*}
		|T_M([\bx],[\by])-T_{M'}([g(\bx)],[g(\by)])|\le 2+2+4=8.
	\end{equation*}
	which implies that
	$\xi^8\rho_M(\bx,\by)\le\rho_{M'}(g(\bx),g(\by))\le\xi^{-8}\rho_M(\bx,\by).$
	Hence $g:\widetilde{\Omega} \to \widetilde{\Omega}$ is a bi-Lipschitz map.

By Lemma \ref{Omegadense},  $\widetilde{\Omega}$
	is dense in ${\mathcal A}_M$, so
	${\mathcal A}_M\simeq {\mathcal A}_{M'}$ by Lemma \ref{extendLip}.
\end{proof}
\section{\textbf{Proof of  Theorem \ref{thm:main}}}\label{sec6}

Recall that a  cross automaton $M$ of Class 2 admits a one-step simplification, so there exists a sequence
$$
M=M_0,\ M_1,\ \dots,\ M_q=M^*
$$
such that $M_j$, $0\leq j<q$ are of Class 2, $M_{j+1}$ is a one-step simplification of $M_j$ and   $M_q$ is of Class $0$. We call  $M^*$  the \emph{final-simplification} of $M$.

By using  Theorem \ref{thm:onesimp} repeatedly, we obtain
\begin{cor}\label{finalspaceLip}
	Let $M$ be a cross automaton of Class 2 and let $M^*$ be the final-simplification of $M$. Then $(\mathcal{A}_M,\rho_M)\simeq(\mathcal{A}_{M^*},\rho_{M^*})$.
\end{cor}

Let $K$ be a Bara\'nski carpet satisfying the cross intersection condition.
If $K$ satisfies the top-isolated condition but not satisfies the vertical separation condition,
we set $M^*_K$ be the final-simplification of $M_K$;
if   $K$  satisfies the vertical separation condition,
we set $M^*_K=M_K$.

\begin{lemma}\label{finalEFLip}
	Let $E$ and $F$ be two Bara\'nski carpets satisfying the assumptions in Theorem \ref{thm:main}.
	Then there exists a map $f:(\mathcal{A}_{M^*_E,\xi},\rho_{M^*_E,\xi})\to(\mathcal{A}_{M^*_F,\xi},\rho_{M^*_F,\xi})$ which is an isometry.
\end{lemma}
\begin{proof}
	Let $I=\{a_1,a_2,\dots,a_k\}\subset\Sigma$ be a H-block of $E$.
	By our assumption in Theorem \ref{thm:main}, there is a size-preserving bijection from the collection of H-blocks of $E$ to that of $F$, which we denote by $\widehat{h}$.
	That is, for any H-block $I=\{a_1,a_2,\dots,a_k\}$ of $E$,
	$$
		\widehat{h}(I)=\{b_1,b_2,\dots,b_k\}
	$$
	is an H-block of $F$, and if $(I_1,I_2)$ is an H-block pair of $E$, then $(\widehat{h}(I_1),\widehat{h}(I_2))$ is an H-block pair of $F$.
	
	Define $h:\Sigma\to\Sigma$ by $h(a_j)=b_j$, that is, if $a_j$ is the $j$-th element of an H-block $I$ of $E$, then define $h(a_j)$ be the $j$-th element of $\widehat{h}(I)$. Clearly, for any $r,s\in\Sigma$,
	\begin{equation}\label{finalEFLip2}
		r\lhd_{H}s \Leftrightarrow h(r)\lhd_{H}h(s)
		\quad{\rm and} \quad
		r\lhd_{\be_1}s\Leftrightarrow  h(r)\lhd_{\be_1}h(s).
	\end{equation}

	Now we define $f:\Sigma^\infty\to\Sigma^\infty$ by
	$
	f\big((x_i)_{i\geq1}\big)=\big(h(x_i)\big)_{i\geq1}.
	$
	Clearly, $f$ is a bijection. Moreover, by \eqref{finalEFLip2}, for any $\bx,\by\in\Sigma^\infty$ we have
	\begin{equation*}
		T_{M_E^*}(\bx,\by)=T_{M_F^*}\big(f(\bx),f(\by)\big).
	\end{equation*}
	It follows that  $[\bx]\mapsto [f(\bx)]$ is an isometry from  $\mathcal{A}_{M^*_E}$ to $\mathcal{A}_{M^*_F}$.
\end{proof}

\begin{proof}[\textbf{Proof of Theorem} \ref{thm:main}]
	By Theorem \ref{thm: Holder map} and Corollary \ref{finalspaceLip}, we have
	\begin{equation*}\label{Lipeq1}
		(E,|\cdot|)\overset{\text{H\"older}}{\simeq} (\mathcal{A}_{M_E,\xi},\rho_{M_E,\xi})\simeq (\mathcal{A}_{M^*_E,\xi},\rho_{M^*_E,\xi})
	\end{equation*}
	and
	\begin{equation*}
		(F,|\cdot|)\overset{\text{H\"older}}{\simeq} (\mathcal{A}_{M_F,\xi},\rho_{M_F,\xi})\simeq (\mathcal{A}_{M^*_F,\xi},\rho_{M^*_F,\xi}),
	\end{equation*}
	Therefore, by Lemma \ref{finalEFLip}, we have $E\overset{\text{H\"older}}{\simeq} F$.
	
	If both $E$ and $F$ are fractal squares, by setting $\xi=1/n$, we obtain that $E\simeq F$.
\end{proof}

\section{\textbf{A universal map}}\label{map g}
Let $N\geq 3$ be an integer. Let $\Sigma=\{1,\dots, N\}$ and let $\gamma, \lambda, \kappa\in \Sigma$ be three distinct letters.
Let $\tau \in \Sigma\setminus \{\gamma, \kappa\}$, that is, it can happen that $\tau=\lambda$.

Set
$$\Omega=\{{\boldsymbol \omega} \kappa^{\infty},{\boldsymbol \omega}\in \Sigma^{*}\}$$
as in Section \ref{one Sim}.
In this section, we recall a bijection map $g:\Omega\to\Omega$ introduced by Huang \emph{et al.} \cite{HWYZ23}.
We emphasize that  the discussion of this section is purely symbolic, and is irrelevant to metric or automaton.


\medskip
\subsection{Segment decomposition}
~\\
\indent  The article \cite{HWYZ23} introduced two decompositions of sequences in $\Sigma$.
Set
\begin{equation}\label{CM}
	\mathcal{C}_M:=\{\tau\gamma^k;k\ge 2\}\cup\{\kappa\lambda^{k}\kappa\gamma;k\ge 0\}.
\end{equation}

\begin{defi}[$M$-initial decomposition]\label{M-segment}
	\emph{Let $\bx=(x_i)_{i=1}^\infty\in\Omega$.  The longest prefix  $X_1$ of $\bx$ satisfying
		$X_1\in {\mathcal C}_M\cup \Sigma$ is called the $M$-initial segment of $\bx$.}
	
	\emph{Inductively, each $\bx=(x_i)_{i=1}^{\infty}\in\Omega$ can be uniquely written as
		$\bx=\prod_{j=1}^\infty X_j:=X_1X_2\cdots X_k \cdots,$
		where $X_k$ is the $M$-initial segment of $\prod_{j\geq k} X_j$. We call $(X_j)_{j\geq 1}$ the \emph{$M$-decomposition} of $\bx$.}
\end{defi}

Next we define $M'$-initial segment and $M'$-decomposition. Set
\begin{equation}\label{CM'}
	\mathcal{C}_{M'}=\{\kappa\lambda^k\kappa\gamma;k\ge 0\}\cup\{\kappa\lambda^{k}\kappa\gamma\gamma;k\ge 0\}\cup\{\tau\gamma\gamma\},
\end{equation}

\begin{defi}[$M'$-initial decomposition]\label{M'-segment}
	\emph{Let $\bu=(u_i)_{i=1}^\infty\in\Omega$.  A word $U_1$  is called the $M'$-initial segment of $\bu$, if
		it is the longest prefix of $\bu$ such that $U_1\in {\mathcal C}_{M'}\cup \Sigma$. Similar as above,   we define the
		\emph{$M'$-decomposition} of $\bu$.
	}
\end{defi}

Two words are said to be \emph{comparable}, if one word is a prefix of another one.

\begin{remark}\label{rem:observe} \emph{Here are two useful observations.}
	
	\emph{ (i) If two elements in ${\mathcal C}_M$ are  comparable,
		then both of them are of the form $\tau\gamma^k$. If two elements in ${\mathcal C}_{M'}$ are  comparable,
		then one of them is $\kappa\lambda^k\kappa\gamma$ and another one is $\kappa\lambda^k\kappa\gamma\gamma$.
	}
	
	\emph{(ii) Let $W\in {\mathcal C}_M\cup {\mathcal C}_{M'}$. Then $W$ is initialled by a word in
		$\{\kappa \lambda, \kappa\kappa, \tau\gamma\}$. Moreover, these words cannot  appear in  $W$ except as a prefix.}
\end{remark}

\medskip
\subsection{Construction of $g$}
~\\
\indent  
First we define $g_0:\mathcal{C}_M\cup\Sigma\rightarrow\mathcal{C}_{M'}\cup\Sigma$ by
\begin{equation}\label{g_0}
g_0:\left\{\begin{array}{rl}
	\tau \gamma^k&\mapsto\kappa \lambda^{k-2}\kappa \gamma,\ k\ge 2;\\
	\kappa \lambda^k\kappa \gamma&\mapsto\kappa\lambda^{k-1}\kappa \gamma\gamma,\ k\ge 1;\\
	\kappa\kappa \gamma&\mapsto\tau \gamma\gamma;\\
	i&\mapsto i,\ \forall~  i\in\Sigma.
\end{array}
\right.
\end{equation}
It is easy to see that $g_0:\mathcal{C}_M\cup\Sigma\rightarrow\mathcal{C}_{M'}\cup\Sigma$ is a bijection. Now we define $g:\Omega\rightarrow\Sigma^\infty$ by
\begin{equation}\label{gxdecomposition}
	g(\bx)=\prod_{j=1}^\infty g_0(X_j),
\end{equation}
where $(X_j)_{j=1}^\infty$ is the $M$-decomposition of $\bx$.

Any $\omega\in \cal C_M$ is not ended with $\kappa$, so  the $M$-decomposition of $\bx=x_1\dots x_k\kappa^\infty\in \Omega$ is $(X_j)_{j=1}^\ell(\kappa)^\infty$  for some $\ell$.
Consequently, $g(\bx)=(\prod_{j=1}^\ell g_0(X_j))(\kappa)^\infty\in\Omega$. Thus $g(\Omega)\subset\Omega$.

\begin{pro}[\cite{HWYZ23}]\label{coincide}
	Let $\bx=x_1x_2\dots, \bu=u_1u_2\dots=g(\bx)$.
	
	\emph{(i)} If $(X_j)_{j=1}^\infty$ is the $M$-decomposition of $\bx$, then the $M'$-decomposition of $g(\bx)$ is $\prod_{j=1}^\infty g_0(X_j)$.
	
	\emph{(ii)} Similarly, if $(U_j)_{j=1}^\infty$ is the $M'$-decomposition of $\bu$, then the $M$-decomposition of $h(\bu)$ is $\left ( g^{-1}_0(U_j) \right )_{j\geq 1}$, where $h(\bu)=\prod_{j=1}^\infty g_0^{-1}(U_j)$.
	
	\emph{(iii)} The map $g:\Omega\rightarrow\Omega$ is a bijection.
\end{pro}

  For $\bx, \by\in\Sigma^\infty$, denote by $\bx\wedge\by$ the maximal common prefix of $\bx$ and $\by\in\Sigma^{\infty}$. For $I\in\Sigma^*$, we denote by $|I|$ the length of $I$. 

\begin{lemma}[\cite{HWYZ23}]\label{lem-length}
	Let $\bx=(x_k)_{k\geq 1}, \by=(y_k)_{k\geq 1}\in\Omega$. Then
	$$|g(\bx)\wedge g(\by)|\geq |\bx\wedge\by|-2.$$
	In other words,  $u_1\cdots u_k$ is determined by $x_1\cdots x_{k+2}$, where $k\ge 1$.
\end{lemma}

The following lemma is our new observation.

\begin{lemma}\label{lem:gamma-free}
	Let $\ba=(a_i)_{i\geq1}\in\Omega$. If $k\geq3$ and $\gamma$ does not occur in $a_1\dots a_k$, then
	$g(\ba)|_{k-2}=\ba|_{k-2}$.
	
\end{lemma}
\begin{proof} Let $(A_i)_{i\geq1}$ be the $M$-decomposition of $\ba$.
	Let $p$ be the largest integer such that    $|A_i|=1$ for all $i\leq p$. (Here $p$ may equal $0$.)
	If $p\geq k-2$, then the lemma holds.
	
	Now we assume that $p\leq k-3$.
	Since $A_{p+1}\in {\mathcal C}_M$, it does not end at or before $a_k$ by the assumption, we have $A_{p+1}\notin\{\tau\gamma^{\ell};\ell\geq2\}\cup\{\kappa\kappa\gamma\}$. So the only choice is
	$A_{p+1}=\kappa\lambda^{\ell}\kappa\gamma$, and we further have  $\ell\geq k-p-2$.  Therefore,
	$$g(\ba)= (a_1\dots a_p) g(A_{p+1})\dots  = (a_1\dots a_p)\kappa\lambda^{\ell-1}\kappa\gamma\gamma\dots,$$
	it follows that 	$|g(\ba)\wedge \ba| \geq p+1+\ell-1\geq k-2$.
	The lemma is proved.
\end{proof}

\section{\textbf{Proof of Theorem \ref{main}}}\label{pf main}

Let $M$ be a  cross automaton of \emph{Class 2} with top and bottom vertices $\gamma$ and $\lambda$ respectively,  let  $M'$ be a one-step simplification of $M$
by deleting $(\tau, \kappa)$, see \eqref{fig_break1}.

For
$\bx, \by\in\Omega=\{\omega \kappa^{\infty};\omega\in\Sigma^*\}$, denote $\bu=g(\bx)~{\rm and}~ \bv=g(\by)$. Let $(X_j)_{j=1}^{\infty}$, $(Y_j)_{j=1}^{\infty}$ be the
$M$-decompositions of $\bx, \by$ respectively, and $(U_j)_{j=1}^{\infty}$, $(V_j)_{j=1}^{\infty}$ be the $M'$-decompositions of $\bu, \bv$ respectively.

In $M$ (or $M'$), for $S,S'\in Q$ and $(i,j)\in\Sigma^2$, we will use $S\stackrel{(i,j)}\longrightarrow S'$ as an alternative notation for $\delta(S,(i,j))=S'$. One should keep in mind that, in both $M$ and $M'$,
\begin{equation}\label{eq:e_2_loop}
	\be_2\stackrel{(i,j)}\longrightarrow \be_2 \textit{\rm ~if and only if ~}(i,j)=(\gamma,\lambda),
\end{equation}
and
  except that  $\be_2\stackrel{(\gamma,\lambda)}\longrightarrow \be_2$, $-\be_2\stackrel{(\lambda, \gamma)}\longrightarrow -\be_2$ and $Id \stackrel{(\gamma,\gamma)}\longrightarrow Id$, we have
\begin{equation}\label{eq:gamma exit}
	 S\stackrel{(\gamma,\theta)}\longrightarrow Exit,\quad
S\stackrel{(\theta,\gamma)}\longrightarrow Exit, \quad \forall \theta\in\Sigma.
\end{equation}
  (See Lemma \ref{def:gamma-isolated}).

\begin{lemma}\label{gotoExit}
	Let $\ba=(a_i)_{i=1}^{\infty},\bb=(b_i)_{i=1}^{\infty}\in\Omega$ with $a_1\neq b_1$.  If $\ba=\lambda^k\kappa\gamma\cdots (k\geq 0)$,
	then
	$$T_M(\ba,\bb), T_{M'}(\ba,\bb)\in \{0,1,2\}.$$
\end{lemma}

\begin{proof}
	Since $M'$ is a simplification of $M$, we have $T_{M'}(\ba,\bb)\leq T_M(\ba,\bb)$.
	So we only need to prove  $T_M(\ba,\bb)\leq2$.
	Suppose $T_M(\ba,\bb)>0$,  then $S_1=\delta(id, (a_1, b_1))\in \{\pm\be_1,\pm\be_2\}$ since $a_1\neq b_1$.

	\ding{172} If $k=0$, then $a_1=\kappa$ and $a_2=\gamma$.
	Since $\kappa$ is $V$-maximal in $M$,
we have $S_1\neq \be_2$, which forces  $S_1=\pm \be_1$. So   $(a_2,b_2)$ leads the state $S_1$ to $Exit$ in $M$
by \eqref{eq:gamma exit}, hence $T_M(\ba,\bb)\leq1$.
	
	\ding{173} If $k=1$, then $a_1=\lambda,a_2=\kappa$ and $a_3=\gamma$.
	Suppose $S=\be_2$, then $(a_2,b_2)=(\gamma,\lambda)$, which contradicts to $a_2=\kappa$.
 So $S_1=\pm \be_1$, and the next state is also $\pm \be_1$ by the self-looping property.
  Thus $\pm \be_1\stackrel{(a_3,b_3)}\longrightarrow Exit$ since $a_3=\gamma$, so $T_M(\ba,\bb)\leq2$.

%
%
	
	\ding{174}  If $k\geq2$, then $a_1=a_2=\lambda$.
	By item (iii) of Lemma \ref{def:gamma-isolated},
   we have
	$id\stackrel{(\lambda\lambda, b_1b_2)}\longrightarrow Exit$, so
  $T_M(\ba,\bb)\leq 1$.
  The lemma is proved.
\end{proof}


Recall that $\mathcal{C}_M=\{\tau \gamma^k;k\ge 2\}\cup\{\kappa \lambda^{k}\kappa \gamma;k\ge 0\}$, see \eqref{CM}.

\begin{lemma}\label{case-id}
	\emph{(i)} Let $\bx,\by\in\Omega$. If $x_1=y_1$ and $X_1\neq Y_1$, then
	\begin{equation}\label{state-id-1}
		T_M(\bx,\by)\leq |\bx\wedge\by|+2.
	\end{equation}
	
	\emph{(ii)} Let $\bu,\bv\in\Omega$. If $u_1= v_1$ and $U_1\neq V_1$, then
	\begin{equation}\label{state-id-2}
		T_{M'}(\bu,\bv)\leq |\bu\wedge\bv|+2.
	\end{equation}
\end{lemma}

\begin{proof} (i)  Let $k=|\bx\wedge\by|$, then $k\ge 1$ and
$$T_M(\bx,\by)=k+T_M(x_{k+1}x_{k+2}\dots,y_{k+1}y_{k+2}\dots).$$
Note that $X_1\neq Y_1$, at least one of $X_1$ and $Y_1$ is in $\mathcal{C}_M$, say $X_1\in\mathcal{C}_M$.
	
\medskip

	\textit{Case 1. }$X_1=\tau \gamma^\ell(\ell\ge 2)$.
	
	In this case, $k\le\ell+1$, for otherwise $X_1=Y_1$, a contradiction. If $k\le\ell$, then $(x_{k+1},y_{k+1})=(\gamma,\theta)$, where $\theta \neq \gamma$; if $k=\ell+1$, then $Y_1=\tau \gamma^s$ with $s\ge \ell+1$ and we have $(x_{k+1},y_{k+1})=(\theta,\gamma)$
 with $\theta\neq \gamma$.  So by \eqref{eq:gamma exit},
 $T_M(x_{k+1}\dots, y_{k+1}\dots)=0$ and    $T_M(\bx,\by)=k$.
	
 	\medskip

	\textit{Case 2. }$X_1=\kappa\lambda^\ell\kappa\gamma(\ell\ge 0)$.
	
	In this case, we have $k\le\ell+2$, for otherwise $X_1=Y_1$.
	If $k\le\ell+1$, then $x_{k+1}x_{k+2}\cdots=\lambda^p\kappa\gamma\cdots (p\geq 0)$,
	by Lemma \ref{gotoExit} we have
$$T_M(\bx,\by)=k+T_M(x_{k+1}x_{k+2}\dots,y_{k+1}y_{k+2}\dots)\in \{k,k+1,k+2\}.$$
	If $k=\ell+2$, then $x_{k+1}=\gamma$ and $(x_{k+1},y_{k+1})=(\gamma,\theta)$ with $\theta\neq \gamma$, so by \eqref{eq:gamma exit} we have $T_M(\bx,\by)=k$.
	Consequently, \eqref{state-id-1} always hlods.
	
	(ii) Using item (i)  we have
	$$
	T_{M'}(\bu,\bv)\leq T_M(\bu,\bv)\leq |\bu\wedge \bv|+2,
	$$
	where the first inequality holds since $M'$ is a simplification of $M$.
\end{proof}

\begin{lemma}\label{lem-key-1}
	Let $\bx,\by\in \Omega$. If $X_1\neq Y_1$, then
	\begin{equation}\label{eq-key-1}
		T_M(\bx,\by)-T_{M'}(\bu,\bv)\le 4.
	\end{equation}
\end{lemma}
\begin{proof}
Let $k=T_M(\bx,\by)$. Obviously, the lemma holds when $k\leq 4$, so in the following we  assume that $k\geq5$.
	Let $S_{M,1}=\delta_M(id, (x_1,y_1))$ be the first state of the itinerary of $(\bx,\by)$ in $M$.
In the following we prove \eqref{eq-key-1} by 3 cases.
	
 	\medskip
	\emph{Case 1}. $S_{M,1}=Id$.\\
	\indent  In this case we have   $x_1=y_1$, so
	$T_M(\bx,\by) \leq |\bx\wedge\by|+2$ by Lemma \ref{case-id} (i).
	Besides, $|\bu\wedge \bv|\geq |\bx\wedge\by|-2$ by Lemma \ref{lem-length}. So
	$$
	T_M(\bx,\by)-T_{M'}(\bu,\bv)\le (|\bx\wedge\by|+2)-|\bu\wedge\bv|\le 4.
	$$
	
 	\medskip
	\emph{Case 2}. $S_{M,1}=\be_1$ or $-\be_1$.\\
	\indent  By symmetry, we   assume that $S_{M,1}=\be_1$.
	The itinerary of $(\bx,\by)$ in $M$  is $id\to (\be_1)^{k}\to Exit$, so
	by  \eqref{eq:gamma exit}, $\gamma$  neither occurs in $x_1\dots x_k$ nor in $y_1\dots y_k$.
	Then by lemma \ref{lem:gamma-free}, we have $u_1\dots u_{k-2}=x_1\dots x_{k-2}$ and $v_1\dots v_{k-2}=y_1\dots y_{k-2}$.
	So $T_{M'}(\bu,\bv)\geq k-2$ since $\cal P'_H=\cal P_H$ and $\cal P'_{\be_1}=\cal P_{\be_1}$,
 and \eqref{eq-key-1} follows.
	
 	\medskip
	
	\emph{Case 3}. $S_{M,1}=\be_2$ or $-\be_2$.\\
	\indent
	By symmetry, we  assume that $S_{M,1}=\be_2$.
	The itinerary of $(\bx,\by)$ in $M$  must be $id\to (\be_2)^{k}\to Exit.$
	Recall that $\cal P_{\be_2}=\{(\gamma,\lambda)\}$,
	so $(x_1, y_1)\in\cal P_V$ and
	\begin{equation*}
		\binom{\bx}{\by}=\binom{x_1\gamma^{k-1}x_{k+1}\dots}{y_1\lambda^{k-1}y_{k+1}\dots}
 	\text{ and }\binom{x_{k+1}}{y_{k+1}}\neq\binom{\gamma}{\lambda}.  
	\end{equation*}
	
	\ding{172} If $x_1\neq \tau$, then $y_1\neq \kappa$ by the uniqueness  property of the cross automaton.
	So $|X_i|=|Y_i|=1$ for all $1\leq i\leq k$, which implies that
	\begin{equation*}
		\binom{\bu}{\bv}=\binom{x_1\gamma^{k-1}\dots}{y_1\lambda^{k-1}\cdots}.
	\end{equation*}
	The itinerary of $(\bu,\bv)$ in $M'$  must be $id\to (\be_2)^{k}\to \cdots$ since $(u_1,v_1)=(x_1,y_1)\in \cal P'_V$.
	Thus $T_{M'}(\bu,\bv)\ge k$ and \eqref{eq-key-1} follows.
	
	\ding{173} If $x_1=\tau$, then $y_1=\kappa$.
	In this case, by the definition of $g$ (see \eqref{gxdecomposition})
	we have
	\begin{equation*}
		\binom{\bu}{\bv}=\binom{\kappa \lambda^{k-3}\dots}{\kappa\lambda^{k-2}\dots}.
	\end{equation*}
	The itinerary of $(\bu,\bv)$ in $M'$  is $id\to (Id)^{k-2}\to \cdots$, so $T_{M'}(\bu,\bv)\ge k-2$ and \eqref{eq-key-1} holds.
\end{proof}

Let $\bu,\bv\in\Omega$. We shall show that  if $U_1\ne V_1$, then
\begin{equation}\label{eq-key-2}
	T_{M'}(\bu,\bv)-T_M(\bx,\by)\leq 4.
\end{equation}
Let $S_{M',1}=\delta_{M'}(id, (u_1,v_1))$ be the first state of the itinerary of $(\bu,\bv)$ in $M'$. If $S_{M',1}=Exit$, it is obvious that \eqref{eq-key-2} holds. We will prove \eqref{eq-key-2} for other cases in the following three lemmas.

Recall that $\mathcal{C}_{M'}=\{\kappa\lambda^k\kappa \gamma;k\ge 0\}\cup\{\kappa\lambda^{k}\kappa \gamma\gamma;k\ge 0\}\cup\{\tau \gamma\gamma\}$, see \eqref{CM'}.

\begin{lemma}\label{lem-key-5}
	Equation \eqref{eq-key-2} holds if $U_1\neq V_1$ and $S_{M',1}=\text{Id}$.
\end{lemma}

\begin{proof}
	That $S_{M',1}=Id$ implies that $u_1=v_1$, so at least one of $U_1$ and $V_1$ is in $\mathcal{C}_{M'}$, say $U_1\in\mathcal{C}_{M'}$.
	
 	\medskip
	\emph{Case 1.} $U_1=\tau \gamma\gamma$.\\
	\indent In this case we have $|\bu\wedge\bv|\le 2$ and $T_{M'}(\bu,\bv)\leq 4$  by Lemma \ref{case-id} (ii).
	So  \eqref{eq-key-2} holds.
	
 	\medskip
	\emph{Case 2.}  $U_1=\kappa\lambda^k\kappa \gamma$ or $\kappa\lambda^k\kappa \gamma\gamma$$(k\geq0)$.\\
	\indent
	\ding{172} If $V_1=\kappa\lambda^\ell\kappa \gamma$ or $\kappa\lambda^\ell\kappa \gamma\gamma(\ell\ge 0)$, it is clear that $\ell\ne k$ when $U_1$ and $V_1$ are of the same type.
	We see that $|\bu\wedge\bv|\le \max\{|U_1|, |V_1|\}\le 4+\min\{k,\ell\}.$
	Applying $g_0^{-1}$ (see \eqref{g_0}) to $U_1$ and $V_1$, we have
	$$
	X_1\in\{\tau \gamma^{k+2},\kappa\lambda^{k+1}\kappa \gamma\},\quad
	Y_1\in\{\tau \gamma^{\ell+2},\kappa\lambda^{\ell+1}\kappa \gamma\}.
	$$
	Then $T_M(\bx,\by)\ge 1+\min\{k+1,\ell+1\}$ holds for all possible combinations of $X_1$ and $Y_1$. Hence
	\begin{equation*}
		T_{M'}(\bu,\bv)-T_M(\bx,\by)\leq |\bu\wedge\bv|+2-T_M(\bx,\by)\leq 4.
	\end{equation*}
	
	\ding{173}
	$V_1=\kappa$, write $\bv$ as $\kappa\lambda^\ell\widetilde{\lambda}\dots$, where $\widetilde{\lambda} \neq \lambda$. Then
	$|\bu\wedge\bv|\leq2+\min\{k,\ell\},$ and $|V_i|=1$ for all $1\leq i\leq \ell$ since $V_1=\kappa$,
 and $\lambda\lambda$ is not a prefix of any words in ${\mathcal C}_{M'}$ (even
 if   $\lambda=\tau$). By the definition of  $g$, see \eqref{gxdecomposition}, we have
	$$
	X_1\in\{\tau \gamma^{k+2},\kappa\lambda^{k+1}\kappa\gamma\},\quad \by=\kappa\lambda^{\ell-1}\cdots.
	$$
	Then $T_M(\bx,\by)\ge 1+\min\{k+1,\ell-1\}$ holds for all possible combinations of $X_1$ and $\by$. Thus
	\begin{equation*}
		T_{M'}(\bu,\bv)-T_M(\bx,\by)\leq |\bu\wedge\bv|+2-T_M(\bx,\by)\leq4.
	\end{equation*}
\end{proof}

\begin{lemma}\label{lem-key-8}
	Equation \eqref{eq-key-2} holds if $S_{M',1}=\be_1$ or $-\be_1$.
\end{lemma}
\begin{proof}
	By symmetry, we   assume that $S_{M',1}=\be_1$.
	Let $k=T_{M'}(\bu,\bv)\geq5$, then
	the itinerary of $(\bu,\bv)$ in $M'$  is
 $id\to (\be_1)^{k}\to Exit.$
 By \eqref{eq:gamma exit}, $\gamma$  neither occurs in $u_1\dots u_k$ nor in $v_1\dots v_k$.

	\medskip

	\emph{Case 1.} $\lambda$ is $\be_1$-isolated.
	
	In this case, $u_i,v_i\notin\{\lambda,\gamma\}(2\leq i\leq k)$, then
	$|U_i|=|V_i|=1$ for all $1\leq i\leq k-2$, hence
	\begin{equation*}
		\binom{\bx}{\by}=\binom{u_1\cdots u_{k-2}x_{k-1}\cdots}{v_1\cdots v_{k-2}y_{k-1}\cdots}.
	\end{equation*}
	The itinerary of $(\bx,\by)$ in $M$  is $id\to (\be_1)^{k-2}\to \cdots$. So $T_M(\bx,\by)\ge k-2$ and \eqref{eq-key-2} holds.
	
	\medskip
	\emph{Case 2.} $\lambda$ is not $\be_1$-isolated.
	
	In this case, there exists a letter $\theta\in \Sigma\setminus\{\gamma\}$ such that $(\lambda,\theta)\in \cal P_{\be_1}\cup \overline{\cal P}_{\be_1}$.

We claim that $\kappa$ is $H$-isolated. Suppose on the contrary that $(\kappa, \eta)\in {\mathcal P}_H$. 	Set
	\begin{equation*}
		\bx'=\kappa\lambda^{\infty},\by'=\eta\theta^{\infty},\bz'=\tau\gamma^{\infty}.
	\end{equation*}
	Then the itineraries of $(\bx',\by')$ and $(\bx',\bz')$ in $M$ are
	$id\to (\be_1)^{\infty}$ and $id\to (-\be_2)^{\infty}$, respectively. So
	$T_M(\bx',\by')=T_M(\bx',\bz')=+\infty$, which contradicts to triple-coding-free condition.
Our claim is proved.

If $|U_1|>1$, then $U_1\in\{\kappa\lambda^{\ell}\kappa\gamma;\ell\geq2\}\cup \{\kappa\lambda^{\ell}\kappa\gamma\gamma;\ell\geq2\}$ since $u_i\neq\gamma$ for all $1\leq i\leq 5$,
so $(\kappa,v_1)=(u_1,v_1)\in {\cal P}_H$,  which contradicts that $\kappa$ is $H$-isolated.
So  $|U_1|=1$.
 By symmetry, we have $|V_1|=1$.
	
 Let $p$ be the largest integer such that   $|U_i|=1$ for all $i\leq p$.
 We claim that $p\geq k-2$.
		Suppose on the contrary  $p\leq k-3$. Then  $U_{p+1}\in\{\kappa\lambda^{\ell}\kappa\gamma;\ell\geq1\}\cup \{\kappa\lambda^{\ell}\kappa\gamma;\ell\geq1\}$ since $u_i\neq \gamma$ for all $p+1\leq i\leq k$. 	
We set
	\begin{equation*}
		\bx''=u_1\dots u_p\kappa\lambda^{\infty},\by''=v_1\dots v_pv_{p+1}\theta^{\infty}
		~\textit{\rm  and }~\bz''=u_1\dots u_p\tau\gamma^{\infty}.
	\end{equation*}
	Obviously, $\bx'',\by''$ and  $\bz''$  are distinct since $u_1\ne v_1$ and  $\theta\ne \gamma$.
	Moreover,  the itineraries of $(\bx'',\by'')$ and $(\bx'',\bz'')$ in $M$ are
	$id\to (\be_1)^{\infty}$ and $id\to(Id)^{p} \to (-\be_2)^{\infty}$ respectively.
	It also contradicts to  triple-coding-free condition.
	
	Therefore,  $x_1\dots x_{k-2}= u_1\dots u_{k-2}$.
By symmetry, we have $y_1\dots y_{k-2}=  v_1\dots v_{k-2}$.
	 By $\cal P'_{H}=\cal P_{H}$ and $\cal P'_{\be_1}=\cal P_{\be_1}$, we obtain $T_{M}(\bx,\by)\geq k-2$, so \eqref{eq-key-2} holds.
\end{proof}

\begin{lemma}\label{lem-key-9}
	Equation \eqref{eq-key-2} holds if $S_{M',1}=\be_2$ or $-\be_2$.
\end{lemma}
\begin{proof}
	By symmetry, we may assume that $S_{M',1}=\be_2$.
	Let $k=T_{M'}(\bu,\bv)\geq5$.
	The itinerary of $(\bu,\bv)$ in $M'$  must be $id\to (\be_2)^{k}\to Exit$.
	So $(u_1,v_1)\in \cal P'_{V}$ and by \eqref{eq:e_2_loop} we have
	\begin{equation*}
		\binom{\bu}{\bv}=\binom{u_1\gamma^{k-1}u_{k+1}\dots}{v_1\lambda^{k-1}v_{k+1}\dots}
		\text{ and }\binom{u_{k+1}}{v_{k+1}}\neq\binom{\gamma}{\lambda}.
	\end{equation*}
 By Lemma \ref{M'carpetauto}, $\kappa$ is $V$-isolated and $\tau$ is
$V$-maximal in $M'$,  so
	  $u_1\notin\{\tau,\kappa\}$ and $v_1\neq \kappa$.
 It follows that  $|U_i|=|V_i|=1$ for all $1\leq i\leq k-1$. Then
 \begin{equation*}
 	\binom{\bx}{\by}=\binom{u_1 \gamma^{k-1}\cdots}{v_1 \lambda^{k-2}\cdots}.
 \end{equation*}
 Since $(u_1,v_1)$ also belongs to ${\mathcal P}_V$, we conclude that $T_M(\bx,\by)\ge k-1$, and \eqref{eq-key-2} follows.
\end{proof}

\begin{proof}[\textbf{Proof of Theorem \ref{main}:}]  Let $\sigma$ be the shift operation  defined by $\sigma((x_k)_{k\geq 1})=(x_k)_{k\geq 2}$.
	For any $\bx,\by\in\Omega$, if $X_1\dots X_k=Y_1\cdots Y_k$ and $X_{k+1}\ne Y_{k+1}$, then $U_1\dots U_k=V_1\cdots V_k$ and $U_{k+1}\ne V_{k+1}$. Let $\ell=|X_1\cdots X_k|$, we have
	$$
		T_M(\bx,\by)-T_{M'}(\bu,\bv)=T_M(\sigma^\ell(\bx),\sigma^\ell(\by))-
		T_{M'}(\sigma^\ell(\bu),\sigma^\ell(\bv)).
	$$
	Since \eqref{eq-key-1} and \eqref{eq-key-2} hold for $\sigma^\ell(\bx),\sigma^\ell(\by),\sigma^\ell(\bu)$ and $\sigma^\ell(\bv)$, we obtain \eqref{eq-key} and we are done.
\end{proof}

\end{document}